\documentclass[11pt]{article}
\usepackage{amsmath}
\usepackage{amssymb}
\usepackage{amsfonts}
\usepackage{color}
\allowdisplaybreaks

\newtheorem{Theorem}{\textbf{Theorem}}[section]
\newtheorem{Lemma}{\textbf{Lemma}}[section]
\newtheorem{Proposition}{\textbf{Proposition}}[section]
\newtheorem{Corollary}{\textbf{Corollary}}[section]
\newtheorem{Remark}{\textbf{Remark}}[section]
\newtheorem{Example}{\textbf{Example}}[section]
\newtheorem{Definition}{\textbf{Definition}}[section]

\newenvironment{theorem}{\begin{Theorem}$\!\!\!$}{\end{Theorem}}
\newenvironment{lemma}{\begin{Lemma}$\!\!\!$}{\end{Lemma}}

\newenvironment{corollary}{\begin{Corollary}$\!\!\!$}{\end{Corollary}}
\newenvironment{remark}{\begin{Remark}$\!\!\!$}{\end{Remark}}

\numberwithin{equation}{section}


\begin{document}
\title{Thresholds for low regularity solutions to \\
wave equations with structural damping}
\author{Tomonori Fukushima and Ryo Ikehata \\
{\small Department of Mathematics, Graduate School of Education, Hiroshima University} \\ {\small Higashi-Hiroshima 739-8524, Japan} \\
and \\
Hironori Michihisa${}^\ast$\\ 
{\small Department of Mathematics, Graduate School of Science, Hiroshima University} \\
{\small Higashi-Hiroshima 739-8526, Japan}
}

\date{%\today
}

\maketitle

\begin{abstract}
We study the asymptotic behavior of solutions to wave equations with a structural damping term 
\[
u_{tt}-\Delta u+\Delta^2 u_t=0,
\qquad
u(0,x)=u_0(x),
\,\,\,
u_t(0,x)=u_1(x),
\]
in the whole space. 
New thresholds are reported in this paper that indicate which of the diffusion wave property and the non-diffusive structure dominates in low regularity cases. 
We develop to that end the previous author's research \cite{FIM} where they have proposed a threshold that expresses whether the parabolic-like property or the wave-like property strongly appears in the solution to some regularity-loss type dissipative wave equation. 
\end{abstract}

\footnote[0]{\hspace{-2em} ${}^\ast$Corresponding author.}
\footnote[0]{\hspace{-2em} \textit{Email:} hi.michihisa@gmail.com}
\footnote[0]{\hspace{-2em} 2010 \textit{Mathematics Subject Classification}. 35B05, 35B33, 35B40, 35B65, 35L30}
\footnote[0]{\hspace{-2em} \textit{Keywords and Phrases}: Structural damping, Regularity-loss, Low regularity, Asymptotic profile, Diffusion wave property, Non-diffusive structure, Threshold}

%\textcolor{red}{赤文字は要修正箇所です. }
%%%%%%%%%%%%%%%%%%%%%%%%
\section{Introduction}
%\label{sec:1}
%%%%%%%%%%%%%%%%%%%%%%%%
In this paper, we study the Cauchy problem of the following wave equation with the structural damping term
\begin{equation}
\label{1.1}
\begin{cases}
u_{tt}-\Delta u+\Delta^2 u_t=0, & t>0,\quad x\in\textbf{R}^n, \\
u(0,x)=u_0(x),
\quad
u_t(0,x)=u_1(x),& x\in\textbf{R}^n, 
\end{cases}
\end{equation}
where $n\ge1$. 
This equation was proposed in \cite{GGH} and they proved that equation~\eqref{1.1} admits a unique mild solution in the class 
\[
u\in C([0,\infty); H^1(\textbf{R}^n))
\cap C^1([0,\infty); L^2(\textbf{R}^n))
\]
if the initial data belong to the energy space 
\begin{align}
\label{1.2}
[u_0,u_1]\in H^1(\textbf{R}^n)\times L^2(\textbf{R}^n).
\end{align}

Before we investigate equation~\eqref{1.1}, we recall previous studies on wave equations with damping terms 
\begin{equation}
\label{1.3}
\begin{cases}
u_{tt}-\Delta u+(-\Delta)^\theta u_t=0, & t>0,\quad x\in\textbf{R}^n, \\
u(0,x)=u_0(x),
\quad
u_t(0,x)=u_1(x),& x\in\textbf{R}^n, 
\end{cases}
\end{equation}
with $\theta\ge0$. 
In the case of $\theta=0$, the classical equation~\eqref{1.3} is known as the damped wave equation. 
In the asymptotic sense, the solution to \eqref{1.3} behaves like a heat kernel. 
See, e.g., D'Abbicco-Reissig \cite{DR}, Han-Milani~\cite{HM}, Hosono~\cite{H}, Hosono-Ogawa~\cite{HO}, Ikehata~\cite{I-0}, Karch~\cite{K}, Marcati-Nishihara~\cite{MN}, Matsumura~\cite{M}, Michihisa~\cite{Mi}, Narazaki~\cite{Na}, Nishihara~\cite{Ni}, Sakata-Wakasugi~\cite{SW} and Takeda~\cite{T}. 
The strongly damped wave equation, i.e., equation~\eqref{1.3} with $\theta=1$, was studied by Ponce \cite{P} and Shibata \cite{S} in the earlier time. 
Results on the asymptotic behavior of the solution can be found in Ikehata~\cite{I}, Ikehata-Onodera~\cite{IO}, Ikehata-Natsume~\cite{IN}, Ikehata-Takeda~\cite{IT}, Ikehata-Todorova-Yordanov \cite{ITY} and  Michihisa~\cite{Mi2, Mi4}. 
Roughly speaking, the solution behaves like the convolution of the heat kernel and the solution to the corresponding wave equation. 
This is so-called a diffusion wave property which is also seen in the case of $\theta\in[1/2,1]$ but not in the lower order damping case $\theta\in[0,1/2)$. 
From these observations, we can deduce that larger values of $\theta$ give stronger wave properties to solutions of equation~\eqref{1.3}.
However, what is common to these cases $\theta\in[0,1]$ is that the norm of Fourier transformed solutions in the high-frequency region is exponentially small. 
That is, only analysis in the low-frequency region is necessary for this concern, and it can be said that the diffusion structure is dominant. 

If $\theta$ is even larger, i.e., $\theta>1$, equation~\eqref{1.3} is of regularity-loss type. 
As a typical case, we are dealing with \eqref{1.1} proposed by Ghisi-Gobbino-Haraux \cite{GGH}. 
On this model, Ikehata-Iyota~\cite{II} derived asymptotic profiles of the solution to \eqref{1.1} with some weighted $L^1$ initial data.
When we consider sufficiently smooth initial data, we can expect the diffusive structure is still dominant even in the case of $\theta>1$. 
Difficulties arise in the low regularity case such as \eqref{1.2}, and then analysis in the high-frequency region is also inevitable to understand the asymptotic behavior of the solution. 
This is just because the high-frequency part of the solution can no longer be regarded as an error. 
Related to the topic, Michihisa~\cite{Mi3} gave higher order asymptotic expansions of the solution to some linear Rosenau-type equation. 
There, we can find the function $e^{-t/|\xi|^2}$ included in terms consisting of the profiles. 
With his technique, quite recently, authors~\cite{FIM} have studied another Rosenau-type equation
\begin{equation*}
\begin{cases}
u_{tt}-\Delta u_{tt}+\Delta^2 u-\Delta u+u_t=0, & t>0,\quad x\in\textbf{R}^n, \\
u(0,x)=u_0(x),
\quad
u_t(0,x)=u_1(x),& x\in\textbf{R}^n, 
\end{cases}
\end{equation*}
where $(u_0,u_1)\in H^{l+1}(\textbf{R}^n)\times H^l(\textbf{R}^n)$ with $l\ge2$. 
This is also of regularity-loss type and they found the leading terms of the solution. 
One is the heat kernel appearing in the low-frequency region whose decay order is determined by the spatial dimension $n$, and the other one is the oscillating functions derived from the high-frequency region whose decay order depends on the regularity of the initial data $l$. 
Compared with their different decay orders, they discovered a meaningful  threshold $l^*=n/2-1$ that indicates whether hyperbolicity or parabolicity is stronger. 
It means that the diffusive structure mainly appears when we impose additional regularity assumptions $l>l^*$ and vice versa. 

In this paper, we give another threshold for such superiority. 
Even if the low regularity Cauchy data are given, we can expect to discuss the similar argument as in \cite{FIM} after pulling out slowly decaying profiles from the solution. 
In this process, we first face the necessity for carrying out higher order expansions of the solution in the high-frequency region as presented in \cite{Mi3}. 
\\

This paper is organized as follows. 
In Section~\ref{sec:2}, we prepare some notation which is commonly used. 
After we define auxiliary functions and some profiles in Section~\ref{sec:3}, we state our results in Section~\ref{sec:4}. 
Proofs of theorems in Section~\ref{sec:4} are written in Section~\ref{sec:5}. 
Results in Section~\ref{sec:6} is the crux of this paper, where we define some new thresholds. 
In Section~\ref{sec:7}, we confirm basic estimates widely used in previous studies. 
There, we also put Lemmas~\ref{lem:7.3}-\ref{lem:7.5} to prove the Theorem~\ref{thm:4.3}. 
Related to these estimates, see also \cite{I}, \cite{IO} and \cite{II}. 

%%%%%%%%%%%
\section{Notation}
\label{sec:2}
%%%%%%%%%%%
Here, we introduce some notation. 
\begin{enumerate}
\item
The set of all positive integers is denoted by $\textbf{N}$ and put $\textbf{N}_0:=\textbf{N}\cup\{0\}$. 
\item
The integer part of $0\le r\in\textbf{R}$ is expressed by $[r]$. 
That is, $[r]:=\max\{k\in\textbf{N}_0: k\le r\}$. 
\item 
The surface area of the unit ball in $\textbf{R}^n$ is expressed by $\omega_n$. 
\item
The Fourier transform $\hat{f}(\xi)$ of a function $f(x)$ is 
\[
\hat{f}(\xi)
:=(2\pi)^{-\frac{n}{2}}
\int_{\textbf{R}^n} 
e^{-ix\cdot\xi}f(x)\,dx.
\]
\item
Throughout this paper, $L^p(\textbf{R}^n)$ represents the usual Lebesgue space and we write its norm as $\|\cdot\|_p$. 
\item 
In connection with the above, we also use the Sobolev space $H^p(\textbf{R}^n)$ equipped with the norm 
\[
\|f\|_{H^p}
:=\left(
\int_{\textbf{R}^n}
(1+|\xi|^{2p})|\hat{f}(\xi)|^2\,d\xi
\right)^\frac{1}{2},
\qquad
f\in H^p(\textbf{R}^n).
\]
%%%
\item 
We define the weighted $L^1$ space as follows: 
\[
L^{1,\gamma}(\textbf{R}^n)
:=\left\{
f\in L^1(\textbf{R}^n):
\|f\|_{1,\gamma}
:=\int_{\textbf{R}^n}
(1+|x|)^\gamma |f(x)|\,dx
<\infty
\right\},
\qquad
\gamma\ge0.
\]
\item 
Let $l\ge0$ and $\gamma\ge0$. 
For $H^{l+1}(\textbf{R}^n)\cap L^{1,\gamma}(\textbf{R}^n)$, $u_1\in H^l(\textbf{R}^n)\cap L^{1,\gamma+1}(\textbf{R}^n)$, put 
\[
I^{l,\gamma}(u_0,u_1)
:=\|u_0\|_{H^{l+1}}
+\|u_1\|_{H^l}
+\|u_0\|_{1,\gamma}
+\|u_1\|_{1,\gamma+1}.
\]

\end{enumerate}

%%%%%%%%%%%%%%
\section{Key functions}
\label{sec:3}
%%%%%%%%%%%%%%
First, we confirm the solution formula (see \cite{Mi5}). 
The characteristic equation corresponding to problem \eqref{1.1} is 
\begin{align*}
\lambda^2+|\xi|^4 \lambda+|\xi|^2=0
\end{align*}
and we put its solutions as 
\begin{align*}
\lambda_1
:=\frac{-|\xi|^4+\sqrt{|\xi|^8-4|\xi|^2}}{2},
\qquad
\lambda_2
:=\frac{-|\xi|^4-\sqrt{|\xi|^8-4|\xi|^2}}{2}.
\end{align*}
Hence, the Fourier transformed solution is formally given by 
\begin{align}
\label{3.1}
\hat{u}(t,\xi)
=E_0(t,\xi)\widehat{u_0}(\xi)
+E_1(t,\xi)
\left(
\frac{|\xi|^4}{2}\widehat{u_0}(\xi)
+\widehat{u_1}(\xi)
\right),
\end{align}
where 
\begin{align*}
E_0(t,\xi):=
\begin{cases}
\displaystyle{
e^{-\frac{1}{2}t|\xi|^4}
\cos
\left(
\frac{t\sqrt{4|\xi|^2-|\xi|^8}}{2}
\right)
},
& |\xi|\le\sqrt[3]{2}, \\[22pt]
\displaystyle{
e^{-\frac{1}{2}t|\xi|^4}
\cosh
\left(
\frac{t\sqrt{|\xi|^8-4|\xi|^2}}{2}
\right)
},
& |\xi|\ge\sqrt[3]{2}, 
\end{cases}
\end{align*}
\begin{align*}
E_1(t,\xi):=
\begin{cases}
\displaystyle{
e^{-\frac{1}{2}t|\xi|^4}
\sin
\left(
\frac{t\sqrt{4|\xi|^2-|\xi|^8}}{2}
\right)
\biggr/
\frac{\sqrt{4|\xi|^2-|\xi|^8}}{2}
},
& |\xi|\le\sqrt[3]{2}, \\[22pt]
\displaystyle{
e^{-\frac{1}{2}t|\xi|^4}
\sinh
\left(
\frac{t\sqrt{|\xi|^8-4|\xi|^2}}{2}
\right)
\biggr/
\frac{\sqrt{|\xi|^8-4|\xi|^2}}{2}
},
& |\xi|\ge\sqrt[3]{2}. 
\end{cases}
\end{align*}
We define 
\[
L_0(t,\xi,a)
:=\cos\left(
t|\xi|-t|\xi|^4 
\cdot
\frac{a}{4+2\sqrt{4-a^2}}
\right),
\]
\[
L_1(t,\xi,a)
:=\sin\left(
t|\xi|-t|\xi|^4 
\cdot
\frac{a}{4+2\sqrt{4-a^2}}
\right)
\biggr/
\frac{|\xi|\sqrt{4-a^2}}{2},
\]
\[
H_0(t,\xi,b)
=\frac{1}{2}
\exp\left(
-\frac{t}{|\xi|^2}
\cdot
\frac{2}{1+\sqrt{1-4b}}
\right),
\]
\[
H_1(t,\xi,b)
=\exp\left(
-\frac{t}{|\xi|^2}
\cdot
\frac{2}{1+\sqrt{1-4b}}
\right)
\biggr/
|\xi|^4 
\sqrt{1-4b}.
\]
Note that 
\[
L_0(t,\xi,0)
=\cos(t|\xi|),
\qquad
L_1(t,\xi,0)
=\frac{\sin(t|\xi|)}{|\xi|},
\]
and 
\[
E_j(t,\xi)
=e^{-\frac{1}{2}t|\xi|^4}
L_j(t,\xi,|\xi|^3),
\qquad
j=0,1.
\]
We also see that 
\[
H_0(t,\xi,|\xi|^{-6})
=\frac{1}{2}
e^{\lambda_1 t}. 
\]
%%%
For $k\in\textbf{N}_0$, we put 
\[
\mathcal{L}_j^k(t,\xi)
:=e^{-\frac{1}{2}t|\xi|^4}
\frac{1}{k!}
\frac{\partial^k L_j}{\partial a^k}(t,\xi,0)
\cdot 
|\xi|^{3k},
\qquad
j=0,1,
\]
\[
m[f]^k(\xi)
:=\sum_{|\alpha|=k}\frac{(-1)^{|\alpha|}}{\alpha!}
\left(
\int_{\textbf{R}^n}
x^\alpha f(x)\, dx
\right)
(i\xi)^\alpha,
\qquad
f\in L^{1,k}(\textbf{R}^n).
\]
To state results of higher order asymptotic expansions in the low-frequency region (see Theorems~\ref{thm:4.1} and \ref{thm:4.2}), 
we prepare the following profiles.  
For $j=0,1$, we define 
\begin{align*}
A_j^k(t,\xi)=
\begin{cases}
B_j^0(t,\xi), & k=0, \\[5pt]
A_j^{k-1}(t,\xi)+B_j^k(t,\xi), & k\in\textbf{N},
\end{cases}
\end{align*}
where
\begin{align*}
B_j^k(t,\xi)
:=\sum_{p=0}^{[k/3]}
\mathcal{L}_j^p(t,\xi)
m[u_j]^{k-3p}(\xi), 
\qquad
k\in\textbf{N}_0.
\end{align*}
That is, 
\begin{align*}
A_j^k(t,\xi)
=\sum_{p=0}^{[k/3]}
\left(
\mathcal{L}_j^p(t,\xi)
\sum_{q=0}^{k-3p}
m[u_j]^q(\xi)
\right).
\end{align*}
%%%
Next, we prepare the following functions leading to higher order asymptotic expansions of the solution in the high-frequency region (see Theorem~\ref{thm:4.6}). 
For $k\in\textbf{N}_0$, we define
\[
\mathcal{H}_j^k(t,\xi)
:=\frac{1}{k!}
\frac{\partial^k H_j}{\partial b^k}(t,\xi,0)\cdot |\xi|^{-6k}, 
\qquad
j=0,1.
\]
For $j=0,1$, we define 
\begin{align*}
C^k(t,\xi)=
\begin{cases}
D^0(t,\xi), & k=0, \\[5pt]
C^{k-1}(t,\xi)+D^k(t,\xi), & k\in\textbf{N},
\end{cases}
\end{align*}
where
\begin{align*}
D^k(t,\xi)
:=
\begin{cases}
\displaystyle{
\left(
\mathcal{H}_0^{k/2}(t,\xi)
+\frac{|\xi|^4}{2}\mathcal{H}_1^{k/2}(t,\xi)
\right)\widehat{u_0}
}, 
& k\in2\textbf{N}_0, \\[11pt]
\displaystyle{
\mathcal{H}_1^{(k-1)/2}(t,\xi)\widehat{u_1}
},
& k\in2\textbf{N}_0+1.
\end{cases}
\end{align*}
That is, 
\begin{align*}
C^k(t,\xi)
:=
\begin{cases}
\displaystyle{
\left(
\mathcal{H}_0^0(t,\xi)
+\frac{|\xi|^4}{2}\mathcal{H}_1^0(t,\xi)
\right)\widehat{u_0}
},
& k=0, \\[11pt]
\displaystyle{
\sum_{p=0}^{[k/2]}
\left(
\mathcal{H}_0^p(t,\xi)
+\frac{|\xi|^4}{2}\mathcal{H}_1^p(t,\xi)
\right)\widehat{u_0}
+\sum_{p=0}^{[(k-1)/2]}
\mathcal{H}_1^p(t,\xi)\widehat{u_1}
}, 
& k\in\textbf{N}.
\end{cases}
\end{align*}

%%%%%%%%%%%%%%%
\section{Results}
\label{sec:4}
%%%%%%%%%%%%%%%
Theorems~\ref{thm:4.1} and \ref{thm:4.2} are results of higher order asymptotic expansions of $E_0(t,\xi)\widehat{u_0}$ and $E_1(t,\xi)\widehat{u_1}$ in the low-frequency region, respectively. 
\begin{theorem}
\label{thm:4.1}
Let $n\ge1$ and $u_0\in L^{1,\gamma}(\textbf{R}^n)$ with $\gamma\ge0$. 
Then, it holds 
\begin{align}
\label{4.1}
\left\|
E_0(t,\xi)\widehat{u_0}-A_0^{[\gamma]}(t,\xi)
\right\|_{L^2(|\xi|\le1)}
\le C\|u_0\|_{1,\gamma}(1+t)^{-\frac{n}{8}-\frac{\gamma}{4}}
\end{align}
for $t\ge0$. 
Here, $C>0$ is a constant independent of $t$ and $u_0$. 
Furthermore, it holds 
\begin{align}
\label{4.2}
\lim_{t\to\infty}
t^{\frac{n}{8}+\frac{\gamma}{4}}
\left\|
E_0(t,\xi)\widehat{u_0}-A_0^{[\gamma]}(t,\xi)
\right\|_{L^2(|\xi|\le1)}
=0.
\end{align}
\end{theorem}

\begin{theorem}
\label{thm:4.2}
Let $n\ge1$ and $u_1\in L^{1,\gamma}(\textbf{R}^n)$ with 
\begin{align}
\label{4.3}
\gamma>\frac{1}{2} \quad (n=1), 
\qquad
\gamma>0 \quad (n=2), 
\qquad
\gamma\ge0 \quad (n\ge3). 
\end{align}
Then, it holds
\begin{align}
\label{4.4}
\left\|
E_1(t,\xi)\widehat{u_1}-A_1^{[\gamma]}(t,\xi)
\right\|_{L^2(|\xi|\le1)}
\le C\|u_1\|_{1,\gamma}(1+t)^{-\frac{n}{8}+\frac{1}{4}-\frac{\gamma}{4}}
\end{align}
for $t\ge0$. 
Here, $C>0$ is a constant independent of $t$ and $u_1$. 
Furthermore, it holds 
\begin{align}
\label{4.5}
\lim_{t\to\infty}
t^{\frac{n}{8}-\frac{1}{4}+\frac{\gamma}{4}}
\left\|
E_1(t,\xi)\widehat{u_1}-A_1^{[\gamma]}(t,\xi)
\right\|_{L^2(|\xi|\le1)}
=0.
\end{align}
\end{theorem}

In the following theorems, we can find some optimal estimates in the low-frequency region, where the diffusive structure is strong in the sense that \eqref{4.2} and \eqref{4.5} holds. 

\begin{theorem}
\label{thm:4.3}
Let $n\ge1$ and $u$ be the solution to \eqref{1.1} with $u_0\in H^1(\textbf{R}^n)\cap L^1(\textbf{R}^n)$, $u_1\in L^2(\textbf{R}^n)\cap L^{1,1}(\textbf{R}^n)$.  
Then, it holds 
\begin{align*}
C_1^1\left|
\int_{\textbf{R}^n} u_1(x)\,dx
\right|
\sqrt{t}
& \le \|\hat{u}(t,\xi)\|_{L^2(|\xi|\le1)}
\le C_2^1 
(\|u_0\|_1+\|u_1\|_{1,1}) 
\sqrt{t}, 
& n=1, \\
C_1^2\left|
\int_{\textbf{R}^n} u_1(x)\,dx
\right|
\sqrt{\log t}
& \le \|\hat{u}(t,\xi)\|_{L^2(|\xi|\le1)}
\le C_2^2 
(\|u_0\|_1+\|u_1\|_{1,1}) 
\sqrt{\log t}, 
& n=2, \\
C_1^n\left|
\int_{\textbf{R}^n} u_1(x)\,dx
\right|
t^{-\frac{n}{8}+\frac{1}{4}}
& \le \|\hat{u}(t,\xi)\|_{L^2(|\xi|\le1)}
\le C_2^n 
(\|u_0\|_1+\|u_1\|_{1,1}) 
t^{-\frac{n}{8}+\frac{1}{4}}, 
& n\ge3,
\end{align*}
for sufficiently large $t$. 
Here, $C_1^n>0$ and $C_2^n>0$ $(n\ge1)$ are constants independent of $t$ and the initial data.   
\end{theorem}

\begin{theorem}
\label{thm:4.4}
Let $n\ge1$ and $u$ be the solution to \eqref{1.1} with $u_0\in H^1(\textbf{R}^n)\cap L^1(\textbf{R}^n)$, $u_1\in L^2(\textbf{R}^n)\cap L^{1,1}(\textbf{R}^n)$. 
Then, it holds 
\begin{align*}
& C_1 \sqrt{
\left(
\int_{\textbf{R}^n} u_0(x)\,dx
\right)^2
+\sum_{j=1}^n
\left(
\int_{\textbf{R}^n} x_j u_1(x)\,dx
\right)^2
}t^{-\frac{n}{8}} \\
& \qquad \le \left\|
\hat{u}(t,\xi)
-\left(
\int_{\textbf{R}^n} u_1(x)\,dx
\right)
e^{-\frac{1}{2}t|\xi|^4}
\frac{\sin(t|\xi|)}{|\xi|}
\right\|_{L^2(|\xi|\le1)}
\le C_2 (\|u_0\|_1+\|u_1\|_{1,1})  t^{-\frac{n}{8}}
\end{align*}
for sufficiently large $t$. 
Here, $C_1>0$ and $C_2>0$ are constants independent of $t$ and the initial data. 
\end{theorem}

\begin{theorem}
\label{thm:4.5}
Let $n\ge1$ and $u$ be the solution to \eqref{1.1} with $u_0\in H^1(\textbf{R}^n)\cap L^{1,1}(\textbf{R}^n)$, $u_1\in L^2(\textbf{R}^n)\cap L^{1,2}(\textbf{R}^n)$. 
Then, it holds 
\begin{align*}
& \left\{
C_1
\sum_{j=1}^n 
\left(
\int_{\textbf{R}^n} x_j^2 u_1(x)\,dx
\right)^2
+C_2\sum_{1\le j<k\le n}
\left(
\int_{\textbf{R}^n} x_j^2 u_1(x)\,dx
\right)
\left(
\int_{\textbf{R}^n} x_k^2 u_1(x)\,dx 
\right) 
\right. \\
& \qquad\qquad
\left. 
+2C_2 
\sum_{1\le j<k\le n}
\left(
\int_{\textbf{R}^n} x_j x_k u_1(x)\,dx
\right)^2
+C_3
\sum_{j=1}^n 
\left(
\int_{\textbf{R}^n} x_j u_0(x)\,dx
\right)^2
\right\}^\frac{1}{2}
t^{-\frac{n}{8}-\frac{1}{4}} \\
\le & \left\|
\hat{u}(t,\xi)
-\left(
\int_{\textbf{R}^n} u_1(x)\,dx
\right)
e^{-\frac{1}{2}t|\xi|^4}
\frac{\sin(t|\xi|)}{|\xi|} 
\right.\\
& \quad \left.
+i \sum_{j=1}^n
\left(
\int_{\textbf{R}^n} x_j u_1(x)\,dx
\right)
e^{-\frac{1}{2}t|\xi|^4}
\frac{\sin(t|\xi|)}{|\xi|}
\xi_j
-\left(
\int_{\textbf{R}^n} u_0(x)\,dx
\right)
e^{-\frac{1}{2}t|\xi|^4}
\cos(t|\xi|)
\right\|_{L^2(|\xi|\le1)} \\
\le & C(\|u_0\|_{1,1}+\|u_1\|_{1,2}) 
t^{-\frac{n}{8}-\frac{1}{4}}
\end{align*}
for sufficiently large $t$. 
Here, we put 
\[
C_1
:=\frac{1}{32}
\int_{|\eta|\le1}
\frac{\eta_1^4}{|\eta|^2}
e^{-|\eta|^4} 
\,d\eta,
\quad
C_2
:=\frac{1}{16}
\int_{|\eta|\le1}
\frac{\eta_1^2 \eta_2^2}{|\eta|^2}
e^{-|\eta|^4} 
\,d\eta,
\quad
C_3
:=\frac{1}{8n}
\int_{|\eta|\le1}
|\eta|^2
e^{-|\eta|^4} 
\,d\eta
\]
and terms including $\displaystyle{\sum_{1\le j<k\le n}}$ are read as zero when $n=1$. 
\end{theorem}

\begin{remark}
{\rm 
Due to the choice of constants $C_j$ $(j=1,2,3)$, the lower bound in the theorem never becomes negative. 
If and only if 
\[
\int_{\textbf{R}^n}
x_j x_k u_1(x)\,dx
=0 
\quad
\mbox{and}
\quad
\int_{\textbf{R}^n}
x_j u_0(x)\,dx
=0 
\quad
\mbox{for all $1\le j\le k\le n$,}
\] 
the lower bound vanishes (see the proof in the next section). 
}
\end{remark}

Next, we state results on the asymptotic behavior of the solution in the high-frequency region. 

\begin{theorem}
\label{thm:4.6}
Let $n\ge1$, $l\ge0$ and $u$ be the solution to \eqref{1.1} with $u_0\in H^{l+1}(\textbf{R}^n)$, $u_1\in H^l(\textbf{R}^n)$. 
For each $k \in \textbf{N}_0$, it holds 
\begin{align}
\label{4.6}
\|\hat{u}(t,\xi)-C^k(t,\xi)\|_{L^2(|\xi| \ge \sqrt{2})}
\le C(\|u_0\|_{H^{l+1}}+\|u_1\|_{H^l})\,(1+t)^{-\frac{l+3k+4}{2}}, \qquad t \ge 0. 
\end{align}
Here, $C>0$ is a constant independent of $t$ and the initial data. 
\end{theorem}

\begin{remark}
{\rm
\begin{enumerate}
\item
When we consider the case of $l=k=0$, 
inequality~\eqref{4.6} becomes 
\begin{align*}
\left\|
\hat{u}(t,\xi)-e^{-\frac{t}{|\xi|^2}}
\widehat{u_0}
\right\|_{L^2(|\xi| \ge \sqrt{2})}
\le C(\|u_0\|_{H^1}+\|u_1\|_2)\,(1+t)^{-2}, \qquad t \ge 0. 
\end{align*}
The profile $e^{-t/|\xi|^2}\widehat{u_0}$ reflects the effect of regularity-loss. 
It decays like $O(t^{-1/2})$ as $t\to\infty$ and thus the above estimate implies that this leading term of the energy solution decays quite slowly. 
\item
We can assure 
\[ 
D^k(t,\xi) \not\equiv 0
\]
for any $k \in \textbf{N}_0$. 
This means that $C^k(t,\xi)\not\equiv C^j(t,\xi)$ if $k\not=j$. 
To confirm the statement, we first consider the case of $k \in 2\textbf{N}_0$. 
Put $k'=k/2$ and 
\[
h(b):=\frac{2}{1+\sqrt{1-4b}}.
\]
The terms in $D^k(t,\xi)$ including the highest order of $t/|\xi|^2$ are
\begin{align*}
& \frac{1}{2 (k')!}
\left(-\frac{t}{|\xi|^2} \right)^{k'}\{h'(0) \}^{k'}\, 
\frac{e^{-\frac{t}{|\xi|^2}}}{|\xi|^{6k'}} \widehat{u_0}(\xi)
+\frac{1}{2 (k')!}
\left(-\frac{t}{|\xi|^2} \right)^{k'}\{h'(0) \}^{k'}\, 
\frac{e^{-\frac{t}{|\xi|^2}}}{|\xi|^{6k'}} \widehat{u_0}(\xi) \\[2mm]
=& \frac{1}{(k')!}
\left(-\frac{t}{|\xi|^2} \right)^{k'}\{h'(0) \}^{k'}\, 
\frac{e^{-\frac{t}{|\xi|^2}}}{|\xi|^{6k'}} \widehat{u_0}(\xi).
\end{align*}
Since $h'(0)=1\not=0$, this value is not identically zero. 
For the same reason, the statement is apparently true for $k\in2\textbf{N}_0+1$. 
\end{enumerate}
}
\end{remark}

The leading term of the solution to \eqref{1.1} is $A_1^0(t,\xi)$ which arises from the low-frequency region or $C^0(t,\xi)$ extracted from the high-frequency region. 
As for the decay orders, the former depends on the spatial dimension $n$ but the latter on the regularity condition $l$. 
The next theorem and the subsequent corollary give a detailed look at the relationship between these two quantities and the asymptotic behavior of the solution to \eqref{1.1}. 
\begin{theorem}
\label{thm:4.7}
Let $n\ge1$, $l\ge0$ and $u$ be the solution to \eqref{1.1} with $u_0\in H^{l+1}(\textbf{R}^n)\cap L^1(\textbf{R}^n)$, $u_1\in H^l(\textbf{R}^n)\cap L^{1,1}(\textbf{R}^n)$. 
Then, it holds 
\begin{align*}
{\rm (i)}\,\,\,& \left\|
\hat{u}(t,\xi)
-e^{-\frac{t}{|\xi|^2}}
\widehat{u_0}
\right\|_2 
\le   
\begin{cases}
CI^{l,0}(u_0,u_1) 
t^{-\frac{l+4}{2}},
& 0\le l\le n/4-9/2 
\,\,\,\mbox{with}\,\,\, n\ge18, \\[7pt]
CI^{l,0}(u_0,u_1) 
t^{-\frac{n}{8}+\frac{1}{4}},
& n/4-9/2< l<n/4-3/2
\,\,\,\mbox{with}\,\,\, n\ge18, \\[7pt]
& \,\,\,\mbox{or}\,\,\,0\le l<n/4-3/2
\,\,\,\mbox{with}\,\,\, 7\le n\le17, 
\end{cases}
\\[11pt]
{\rm (ii)}\,\,\,& \left\|
\hat{u}(t,\xi)
-\left(
\int_{\textbf{R}^n} u_1(x)\,dx
\right)
e^{-\frac{1}{2}t|\xi|^4}
\frac{\sin(t|\xi|)}{|\xi|}
-e^{-\frac{t}{|\xi|^2}}
\widehat{u_0}
\right\|_2 \\[11pt]
& \qquad\qquad\qquad\qquad\quad 
\le CI^{l,0}(u_0,u_1) 
t^{-\frac{n}{8}}, 
\qquad
l=n/4-3/2 \,\,\,\mbox{with}\,\,\, n\ge6,
\\[11pt]
{\rm (iii)}\,\,\,& \left\|
\hat{u}(t,\xi)
-\left(
\int_{\textbf{R}^n} u_1(x)\,dx
\right)
e^{-\frac{1}{2}t|\xi|^4}
\frac{\sin(t|\xi|)}{|\xi|}
\right\|_2 \\[11pt]
& \qquad\qquad\qquad\qquad\quad
 \le 
\begin{cases}
C I^{l,0}(u_0,u_1) t^{-\frac{l+1}{2}}, 
& n/4-3/2<l\le n/4-1 
\,\,\,\mbox{with}\,\,\,n\ge6, \\[7pt]
& \,\,\,\mbox{or}\,\,\,0\le l\le n/4-1 
\,\,\,\mbox{with}\,\,\,n=4, 5, \\[7pt]
C I^{l,0}(u_0,u_1) t^{-\frac{n}{8}}, 
& n/4-1< l 
\,\,\,\mbox{with}\,\,\,n\ge4, \\[7pt]
& \,\,\,\mbox{or}\,\,\,0\le l 
\,\,\,\mbox{with}\,\,\,n=1,2,3, 
\end{cases}
\end{align*}
for $t\ge1$. 
Here, $C>0$ is a constant independent of $t$ and the initial data. 
\end{theorem}

\begin{corollary}
\label{cor:4.1}
\begin{enumerate}
\item Let $n\ge6$ and $0\le l\le n/4-3/2$. 
If $u_0\in H^{l+1}(\textbf{R}^n)\cap L^1(\textbf{R}^n)$, $u_1\in H^l(\textbf{R}^n)\cap L^{1,1}(\textbf{R}^n)$, 
then the solution $u$ to \eqref{1.1} satisfies
\begin{align*}
\|u(t,\cdot)\|_2
\le C I^{l,0}(u_0,u_1)
t^{-\frac{l+1}{2}},
\qquad
t\gg1.
\end{align*}
Here, $C>0$ is a constant independent of $t$ and the initial data. 
\item Let $n\ge1$ and 
\begin{align*}
\begin{cases}
l\ge0, & 1\le n\le 6, \\[7pt]
l\ge n/4-3/2, & n\ge7. 
\end{cases} 
\end{align*}
If $u_0\in H^{l+1}(\textbf{R}^n)\cap L^1(\textbf{R}^n)$, $u_1\in H^l(\textbf{R}^n)\cap L^{1,1}(\textbf{R}^n)$, 
then the solution $u$ to \eqref{1.1} satisfies
\begin{align*}
\|u(t,\cdot)\|_2
\le 
\begin{cases}
C I^{l,0}(u_0,u_1)
\sqrt{t}, 
& n=1, \\[7pt]
C I^{l,0}(u_0,u_1)
\sqrt{\log t}, 
& n=2, \\[7pt]
C I^{l,0}(u_0,u_1)
t^{-\frac{n}{8}+\frac{1}{4}}, 
& n\ge3,
\end{cases}
\qquad
t\gg1.
\end{align*}
Here, $C>0$ is a constant independent of $t$ and the initial data. 
\end{enumerate}
\end{corollary}

\begin{remark}
\label{rem:4.3}
{\rm 
The value $n/4-3/2$ itself was already found in \cite{II}. 
They obtained the optimal estimate for the solution to \eqref{1.1} under $l>n/4-3/2$ but the asymptotic behavior of the solution was not investigated. 
However, from Theorem~\ref{thm:4.7} and Corollary~\ref{cor:4.1}, we can characterize the exponent $l^*:=n/4-3/2$ as a threshold that indicates the superiority of the diffusion wave property or the non-diffusive structute. 
This viewpoint comes from \cite{FIM}, which leads to the results in Section~\ref{sec:6} below. 
}
\end{remark}

%%%%%%%%%%%%%
\section{Proofs}
\label{sec:5}
%\label{sec:}
%%%%%%%%%%%%%
We confirm the following results on expansions of evolution operator $E_j(t,\xi)$ $(j=0,1)$. 
Similar results were already found in \cite{Mi2}. 
\begin{lemma}
\label{lem:5.1}
Let $n\ge1$ and $k\in\textbf{N}_0$. 
Then, there exists a constant $C>0$ such that 
\begin{align}
\label{5.1}
\left\|
E_0(t,\xi)
-\sum_{p=0}^k
\mathcal{L}_0^p(t,\xi)
\right\|_{L^2(|\xi|\le1)}
\le C(1+t)^{-\frac{n}{8}-\frac{3(k+1)}{4}},
\end{align}
\begin{align}
\label{5.2}
\left\|
E_1(t,\xi)
-\sum_{p=0}^k
\mathcal{L}_1^p(t,\xi)
\right\|_{L^2(|\xi|\le1)}
\le C(1+t)^{-\frac{n}{8}+\frac{1}{4}-\frac{3(k+1)}{4}},
\end{align}
for $t\ge0$.
\end{lemma}
\textbf{Proof.} 
We prove \eqref{5.2} only since \eqref{5.1} can be proved in the same way. 
The Taylor theorem gives 
\[
E_1(t,\xi)
-\sum_{p=0}^k
\mathcal{L}_1^p(t,\xi)
=e^{-\frac{1}{2}t|\xi|^4}
\frac{1}{(k+1)!}
\frac{\partial^{k+1} L_1}{\partial a^{k+1}}(t,\xi,\tau|\xi|^3)
\cdot 
|\xi|^{3(k+1)}
\]
for some $0\le\tau\le1$. 
Note that the function $\sqrt{4-a^2}$ and its derivatives are all bounded for $0\le a\le1$.  
Thus, for each $j\in\textbf{N}_0$, there exists a constant $C>0$ such that 
\begin{align}
\label{5.3}
\left|
\frac{\partial^j L_1}{\partial a^j}(t,\xi,a)
\right|
\le \frac{C}{|\xi|}
\sum_{p=0}^j
(t|\xi|^4)^p
\end{align}
for $t\ge0$, $\xi\in\textbf{R}^n$ and $0\le a\le1$. 
Now, we consider $\xi\in\textbf{R}^n$ with $|\xi|\le1$, it follows from \eqref{7.1} that 
\begin{align*}
\left\|
E_1(t,\xi)
-\sum_{p=0}^k
\mathcal{L}_1^p(t,\xi)
\right\|_{L^2(|\xi|\le1)}
& \le C\sum_{p=0}^k
\left\|
(t|\xi|^4)^p
|\xi|^{3(k+1)-1}
e^{-\frac{1}{2}t|\xi|^4}
\right\|_{L^2(|\xi|\le1)} \\
& \le C(1+t)^{-\frac{n}{8}-\frac{3(k+1)-1}{4}}, 
\qquad
t\ge0. 
\end{align*}
The proof of \eqref{5.2} is complete. 
$\hfill \Box$
\\

Based on Lemma~\ref{lem:5.1}, we can obtain \eqref{4.1} and \eqref{4.4}. 
The proof of Theorem~\ref{thm:4.2} is more complicated than that of Theorem~\ref{thm:4.1}. 
So, we give the proof of Theorem~\ref{thm:4.2} briefly. \\
\textbf{Proof of Theorem~\ref{thm:4.2}.} 
First, we rearrange 
\begin{align*}
& E_1(t,\xi)\widehat{u_1}(\xi) \\
= & \left\{
\sum_{p=0}^{[\gamma/3]}
\mathcal{L}_1^p(t,\xi)
+\left(
E_1(t,\xi)
-\sum_{p=0}^{[\gamma/3]}
\mathcal{L}_1^p(t,\xi)
\right)
\right\}
\times
\left\{
\sum_{q=0}^{[\gamma]}
m[u_1]^q(\xi)
+\left(
\widehat{u_1}(\xi)
-\sum_{q=0}^{[\gamma]}
m[u_1]^q(\xi)
\right)
\right\} \\
=& \left(
\sum_{p=0}^{[\gamma/3]}
\mathcal{L}_1^p(t,\xi)
\right)
\left(
\sum_{q=0}^{[\gamma]}
m[u_1]^q(\xi)
\right) \\
& \quad
+\left(
\sum_{p=0}^{[\gamma/3]}
\mathcal{L}_1^p(t,\xi)
\right)
\left(
\widehat{u_1}(\xi)
-\sum_{q=0}^{[\gamma]}
m[u_1]^q(\xi)
\right)
+\left(
E_1(t,\xi)
-\sum_{p=0}^{[\gamma/3]}
\mathcal{L}_1^p(t,\xi)
\right)
\widehat{u_1}(\xi).
\end{align*}
It follows from \eqref{5.3} and \eqref{7.3} with \eqref{7.1} that 
\begin{align*}
\left\|
\left(
\sum_{p=0}^{[\gamma/3]}
\mathcal{L}_1^p(t,\xi)
\right)
\left(
\widehat{u_1}(\xi)
-\sum_{q=0}^{[\gamma]}
m[u_1]^q(\xi)
\right)
\right\|_{L^2(|\xi|\le1)}
& \le C\|u_1\|_{1,\gamma} 
\sum_{p=0}^{[\gamma/3]}
\left\|
(t|\xi|^4)^p
|\xi|^{\gamma-1}
e^{-\frac{1}{2}t|\xi|^4}
\right\|_{L^2(|\xi|\le1)} \\
& \le C \|u_1\|_{1,\gamma}
\left\||\xi|^{\gamma-1}e^{-\frac{1}{4}t|\xi|^4}\right\|_{L^2(|\xi|\le1)} \\
& \le C\|u_1\|_{1,\gamma} (1+t)^{-\frac{n}{8}+\frac{1}{4}-\frac{\gamma}{4}},
\qquad t\ge0. 
\end{align*}
Here, we used \eqref{4.3} to assure the following integrability:
\begin{align*}
\left\||\xi|^{\gamma-1}e^{-\frac{1}{4}t|\xi|^4}\right\|_{L^2(|\xi|\le1)}^2
& \le C(1+t)^{-\frac{n}{4}+\frac{1}{2}-\frac{\gamma}{2}}
\int_{|\eta|\le t^\frac{1}{4}}
|\eta|^{2\gamma-2}e^{-\frac{1}{2}|\eta|^4}\,d\eta \\
& \le C(1+t)^{-\frac{n}{4}+\frac{1}{2}-\frac{\gamma}{2}}
\int_0^\infty s^{n+2\gamma-3}e^{-\frac{1}{2}s^4}\,ds \\
& \le C(1+t)^{-\frac{n}{4}+\frac{1}{2}-\frac{\gamma}{2}}, 
\qquad t\ge0. 
\end{align*}
From \eqref{5.2}, we have 
\begin{align*}
\left\|
\left(
E_1(t,\xi)
-\sum_{p=0}^{[\gamma/3]}
\mathcal{L}_1^p(t,\xi)
\right)
\widehat{u_1}(\xi)
\right\|_{L^2(|\xi|\le1)}
& \le C(1+t)^{-\frac{n}{8}+\frac{1}{4}-\frac{3([\gamma/3]+1)}{4}}
\|u_1\|_1 \\
& \le C(1+t)^{-\frac{n}{8}+\frac{1}{4}-\frac{\gamma}{4}}\|u_1\|_1,
\qquad
t\ge0.
\end{align*}
Later, we deal with the term 
\[
\left(
\sum_{p=0}^{[\gamma/3]}
\mathcal{L}_1^p(t,\xi)
\right)
\left(
\sum_{q=0}^{[\gamma]}
m[u_1]^q(\xi)
\right).
\]
If $[\gamma/3]=0$, then it becomes $A_1^0(t,\xi)$ and \eqref{4.4} is obtained. 
Next, we consider the case of $[\gamma/3]>0$, that is, $\gamma\ge3$. 
Then, we see that 
\[
\left(
\sum_{p=0}^{[\gamma/3]}
\mathcal{L}_1^p(t,\xi)
\right)
\left(
\sum_{q=0}^{[\gamma]}
m[u_1]^q(\xi)
\right)
=A_1^{[\gamma]}(t,\xi)
+\sum_{p=1}^{[\gamma/3]}
\left(
\mathcal{L}_1^p(t,\xi)
\sum_{q=[\gamma]-3p+1}^{[\gamma]}
m[u_1]^q(\xi)
\right)
\]
Here, we note that 
\[
\left[
\frac{[\gamma]}{3}
\right]
=\left[
\frac{\gamma}{3}
\right]
\]
for all $\gamma\ge0$. 
When we confirm 
\begin{align*}
\left\|
\sum_{p=1}^{[\gamma/3]}
\left(
\mathcal{L}_1^p(t,\xi)
\sum_{q=[\gamma]-3p+1}^{[\gamma]}
m[u_1]^q(\xi)
\right)
\right\|_{L^2(|\xi|\le1)} 
& \le C\|u_1\|_{1,[\gamma]}
\sum_{p=1}^{[\gamma/3]}
\left\|
(t|\xi|^4)^p
|\xi|^{[\gamma]}
e^{-\frac{1}{2}t|\xi|^4}
\right\|_{L^2(|\xi|\le1)} \\
& \le C\|u_0\|_{1,[\gamma]}
(1+t)^{-\frac{n}{8}-\frac{[\gamma]}{4}},
\qquad 
t\ge0, 
\end{align*}
we obtain \eqref{4.4}. 

Estimate~\eqref{4.5} can be proved by the Lebesgue dominated convergence theorem (see \cite{IM, Mi4} in detail). 
$\hfill \Box$
\\

Next, we prove Theorem~\ref{thm:4.7}. 
Here, we closely compare two different decay orders of $A_1^0(t,\xi)$ and $C^0(t,\xi)$. 
This motivation and technique were first proposed in \cite{FIM} and recall their threshold $l^*=n/2-1$ stated in our Introduction. 
\\

The proofs of Thereoms~\ref{thm:4.4} and \ref{thm:4.5} are essentially relied on direct calculations. 
See also the corresponding proofs in \cite{Mi4}. 
\\
\textbf{Proof of Theorem~\ref{thm:4.4}.} 
We see that \eqref{4.4} with $\gamma=1$ and the estimate 
\begin{align*}
%\label{5.4}
\left\|
B_1^1(t,\xi)
\right\|_{L^2(|\xi|\le1)}
\le C\|u_1\|_{1,1}(1+t)^{-\frac{n}{8}}, 
\qquad
t\ge0,
\end{align*}
give the upper bound. 
So, we prove the lower bound. 
First, it follows that 
\begin{align*}
& \left\|
\hat{u}(t,\xi)
-\left(
\int_{\textbf{R}^n} u_1(x)\,dx
\right)
e^{-\frac{1}{2}t|\xi|^4}
\frac{\sin(t|\xi|)}{|\xi|}
\right\|_{L^2(|\xi|\le1)}
=\left\|
\hat{u}(t,\xi)
-A_1^0(t,\xi)
\right\|_{L^2(|\xi|\le1)} \\[5pt]
\ge & \left\|
B_1^1(t,\xi)+B_0^0(t,\xi)
\right\|_{L^2(|\xi|\le1)}
-\left\|
E_1(t,\xi)\widehat{u_0}-A_1^1(t,\xi)
\right\|_{L^2(|\xi|\le1)} \\
& \qquad
-\left\|
E_0(t,\xi)\widehat{u_0}-A_0^0(t,\xi)
\right\|_{L^2(|\xi|\le1)} 
-\left\|
\frac{|\xi|^4}{2}E_1(t,\xi)
\widehat{u_0}
\right\|_{L^2(|\xi|\le1)} \\
\ge & \left\|
B_1^1(t,\xi)+B_0^0(t,\xi)
\right\|_{L^2(|\xi|\le1)}
-o(t^{-\frac{n}{8}})
-o(t^{-\frac{n}{8}})
-O(t^{-\frac{n}{8}-\frac{3}{4}}),
\qquad
t\to\infty. 
\end{align*}
Here, we used \eqref{4.5} with $\gamma=1$ and \eqref{4.2} with $\gamma=0$.
Direct calculation shows
\begin{align*}
B_1^1(t,\xi)
=\mathcal{L}_1^0(t,\xi)m[u_1]^1(\xi) 
= -i 
\sum_{j=1}^n
\left(
\int_{\textbf{R}^n} x_j u_1(x)\,dx
\right)
e^{-\frac{1}{2}t|\xi|^4}
\frac{\sin(t|\xi|)}{|\xi|}
\xi_j,
\end{align*}
\begin{align*}
B_0^0(t,\xi)
=\mathcal{L}_0^0(t,\xi)m[u_0]^0(\xi) 
=\left(
\int_{\textbf{R}^n} u_0(x)\,dx
\right)
e^{-\frac{1}{2}t|\xi|^4}
\cos(t|\xi|).
\end{align*}
If $B_1^1(t,\xi)+B_0^0(t,\xi)\equiv0$, then the lower bound in the theorem becomes zero. 
Thus, it suffices to show it under $B_1^1(t,\xi)+B_0^0(t,\xi)\not\equiv0$. 
In this case, one has 
\begin{align*}
& \left\|
B_1^1(t,\xi)+B_0^0(t,\xi) 
\right\|_{L^2(|\xi|\le1)}^2 
=\left\|
B_1^1(t,\xi) 
\right\|_{L^2(|\xi|\le1)}^2
+\left\|
B_0^0(t,\xi) 
\right\|_{L^2(|\xi|\le1)}^2 \\
= & \int_{|\xi|\le1}
e^{-t|\xi|^4} 
\frac{\sin^2(t|\xi|)}{|\xi|^2}
\biggr(
\sum_{j=1}^n 
P_{1,j} \xi_j
\biggr)^2
\,d\xi 
+P_0^2 
\int_{|\xi|\le1}
e^{-t|\xi|^4} 
\cos^2(t|\xi|)
\,d\xi.
\end{align*}
Here, we put 
\[
P_{1,j}
:=\int_{\textbf{R}^n} x_j u_1(x)\,dx, 
\qquad
P_0:=\int_{\textbf{R}^n} u_0(x)\,dx.
\]
Since 
\[
\int_{|\xi|\le1}
e^{-t|\xi|^4} 
\frac{\sin^2(t|\xi|)}{|\xi|^2}
\xi_j \xi_k
\,d\xi
=0
\]
for $1\le j<k\le n$ with $n\ge2$, 
one has 
\begin{align*}
& \int_{|\xi|\le1}
e^{-t|\xi|^4} 
\frac{\sin^2(t|\xi|)}{|\xi|^2}
\biggr(
\sum_{j=1}^n 
P_{1,j} \xi_j
\biggr)^2
\,d\xi 
=\sum_{j=1}^n 
P_{1,j}^2
\int_{|\xi|\le1}
e^{-t|\xi|^4} 
\frac{\sin^2(t|\xi|)}{|\xi|^2}
\xi_j^2
\,d\xi \\
= & \frac{1}{n}
\left(
\int_{|\xi|\le1}
e^{-t|\xi|^4} 
\sin^2(t|\xi|)
\,d\xi
\right)
\sum_{j=1}^n 
P_{1,j}^2 
\ge \frac{1}{n}
\left(
\int_{|\eta|\le1}
e^{-|\eta|^4} 
\sin^2(t^\frac{3}{4}|\eta|)
\,d\eta
\right)
\left(
\sum_{j=1}^n 
P_{1,j}^2 
\right)
t^{-\frac{n}{4}} \\
\ge &  \frac{1}{4n}
\left(
\int_{|\eta|\le1}
e^{-|\eta|^4} 
\,d\eta
\right)
\left(
\sum_{j=1}^n 
P_{1,j}^2
\right)
t^{-\frac{n}{4}},
\qquad t\gg1,
\end{align*}
for all spatial dimensions $n\ge1$. 
Now, we use the following estimate derived by the Riemann-Lebesgue lemma whose origin is in \cite{I}:
\begin{align*}
\int_{|\eta|\le1}
e^{-|\eta|^4} 
\sin^2(t^\frac{3}{4}|\eta|)
\,d\eta
& =\frac{1}{2}
\int_{|\eta|\le1}
e^{-|\eta|^4} 
\,d\eta
-\frac{1}{2}
\int_{|\eta|\le1}
e^{-|\eta|^4} 
\cos(2t^\frac{3}{4}|\eta|)
\,d\eta \\
& \ge \frac{1}{4}
\int_{|\eta|\le1}
e^{-|\eta|^4} 
\,d\eta,
\qquad
t\gg1. 
\end{align*}
Similarly, we have  
\[
P_0^2 
\int_{|\xi|\le1}
e^{-t|\xi|^4} 
\cos^2(t|\xi|)
\,d\xi
\ge 
\frac{1}{4}
\left(
\int_{|\eta|\le1}
e^{-|\eta|^4}
\,d\xi
\right)
P_0^2
t^{-\frac{n}{4}},
\qquad
t\gg1. 
\]
Therefore, we obtain the theorem. 
$\hfill \Box$
\\

\noindent
\textbf{Proof of Theorem~\ref{thm:4.5}.} 
We prove the lower bound only. 
It follows that 
\begin{align*}
& \left\|
\hat{u}(t,\xi)
-\left(
\int_{\textbf{R}^n} u_1(x)\,dx
\right)
e^{-\frac{1}{2}t|\xi|^4}
\frac{\sin(t|\xi|)}{|\xi|} 
\right.\\
& \quad \left.
+i \sum_{j=1}^n
\left(
\int_{\textbf{R}^n} x_j u_1(x)\,dx
\right)
e^{-\frac{1}{2}t|\xi|^4}
\frac{\sin(t|\xi|)}{|\xi|}
\xi_j
-\left(
\int_{\textbf{R}^n} u_0(x)\,dx
\right)
e^{-\frac{1}{2}t|\xi|^4}
\cos(t|\xi|)
\right\|_{L^2(|\xi|\le1)} \\
= &  \left\|
\hat{u}(t,\xi)
-A_1^1(t,\xi)
-A_0^0(t,\xi)
\right\|_{L^2(|\xi|\le1)} \\
\ge & \left\|
B_1^2(t,\xi)+B_0^1(t,\xi)
\right\|_{L^2(|\xi|\le1)} 
-\left\|
E_1(t,\xi)\widehat{u_1}
-A_1^2(t,\xi)
\right\|_{L^2(|\xi|\le1)} \\
& \quad
-\left\|
E_0(t,\xi)\widehat{u_0}
-A_0^1(t,\xi)
\right\|_{L^2(|\xi|\le1)} 
-\left\|
\frac{|\xi|^4}{2}E_1(t,\xi)
\widehat{u_0}
\right\|_{L^2(|\xi|\le1)} \\
\ge & \left\|
B_1^2(t,\xi)+B_0^1(t,\xi)
\right\|_{L^2(|\xi|\le1)}
-o(t^{-\frac{n}{8}-\frac{1}{4}})
-o(t^{-\frac{n}{8}-\frac{1}{4}})
-O(t^{-\frac{n}{8}-\frac{3}{4}}),
\qquad
t\to\infty. 
\end{align*}
Here, we used \eqref{4.5} with $\gamma=2$ and \eqref{4.2} with $\gamma=1$. 
For the same reason as in the previous proof, we consider the case of $B_1^2(t,\xi)+B_0^1(t,\xi)\not\equiv0$. 
By the definition, we have 
\begin{align*}
B_1^2(t,\xi)
=\mathcal{L}_1^0(t,\xi)m[u_1]^2(\xi) 
= e^{-\frac{1}{2}t|\xi|^4}
\frac{\sin(t|\xi|)}{|\xi|}
m[u_1]^2(\xi),
\end{align*}
\begin{align*}
B_0^1(t,\xi)
=\mathcal{L}_0^0(t,\xi)m[u_1]^1(\xi) 
=-i 
\sum_{j=1}^n
\left(
\int_{\textbf{R}^n} x_j u_0(x)\,dx
\right)
e^{-\frac{1}{2}t|\xi|^4}
\cos(t|\xi|)
\xi_j. 
\end{align*}
First, we deal with $B_1^2(t,\xi)$. 
Since $B_1^2(t,c\xi)=c^2 B_1^2(t,\xi)$ for all $c\in\textbf{R}$, we have 
\begin{align*}
\left\|
B_1^2(t,\xi) 
\right\|_{L^2(|\xi|\le1)}^2 
& =\int_{|\xi|\le1}
e^{-t|\xi|^4} \frac{\sin^2(t|\xi|)}{|\xi|^2} 
\big(
m[u_1]^2(\xi)
\big)^2
\, d\xi \\
& \ge t^{-\frac{n}{4}-\frac{1}{2}}
\int_{|\eta|\le1}
e^{-|\eta|^4} \frac{\sin^2(t^\frac{3}{4}|\eta|)}{|\eta|^2} 
\big(
m[u_1]^2(\eta)
\big)^2
\, d\eta,
\qquad
t\ge1.
\end{align*}
If $m[u_1]^2(\eta)\not\equiv0$, the Riemann-Lebesgue lemma gives 
\begin{align*}
\left\|
B_1^2(t,\xi) 
\right\|_{L^2(|\xi|\le1)}^2
\ge \frac{1}{4}t^{-\frac{n}{4}-\frac{1}{2}}
\int_{|\eta|\le1}
e^{-|\eta|^4} \frac{1}{|\eta|^2}
\big(
m[u_1]^2(\eta)
\big)^2
\, d\eta,
\qquad 
t\gg1.
\end{align*}
This estimate, however, holds even if $m[u_1]^2(\eta)\equiv0$. 
When $n=1$, we see that 
\[
m[u_1]^2(\eta)
=-\frac{1}{2}
\left(
\int_{-\infty}^\infty x^2 u_1(x)\,dx
\right) 
\eta^2
\]
and thus 
\[
\left\|
B_1^2(t,\xi) 
\right\|_{L^2(|\xi|\le1)}^2
\ge 2C_1
\left(
\int_{-\infty}^\infty x^2 u_1(x)\,dx
\right)^2
t^{-\frac{1}{4}-\frac{1}{2}},
\qquad
t\gg1.
\]
On the other hand, the case of $n\ge2$ is more complicated. 
Since 
\[
m[u_1]^2(\eta)
=-\frac{1}{2}\sum_{j=1}^n
\left(
\int_{\textbf{R}^n} x_j^2 u_1(x)\,dx
\right)
\eta_j^2 \\
-\sum_{1\le j<k\le n}
\left(
\int_{\textbf{R}^n} x_j x_k u_1(x)\,dx
\right)
\eta_j \eta_k, 
\]
we have 
\begin{align*}
& \int_{|\eta|\le1}
e^{-|\eta|^4} \frac{1}{|\eta|^2}
\big(
m[u_1]^2(\eta)
\big)^2
\, d\eta \\
= & \frac{1}{4}
\sum_{j=1}^n 
\left(
\int_{\textbf{R}^n} x_j^2 u_1(x)\,dx
\right)^2
\int_{|\eta|\le1}
\frac{\eta_j^4}{|\eta|^2}
e^{-|\eta|^4}
\,d\eta \\
& \quad
+\frac{1}{2}
\sum_{1\le j<k\le n}
\left(
\int_{\textbf{R}^n} x_j^2 u_1(x)\,dx
\right)
\left(
\int_{\textbf{R}^n} x_k^2 u_1(x)\,dx
\right)
\int_{|\eta|\le1}
\frac{\eta_j^2 \eta_k^2}{|\eta|^2}
e^{-|\eta|^4} 
\,d\eta \\
& \quad
+\sum_{1\le j<k\le n}
\left(
\int_{\textbf{R}^n} x_j x_k u_1(x)\,dx
\right)^2
\int_{|\eta|\le1}
\frac{\eta_j^2 \eta_k^2}{|\eta|^2}
e^{-|\eta|^4} 
\,d\eta \\
= & 8C_1
\sum_{j=1}^n 
\left(
\int_{\textbf{R}^n} x_j^2 u_1(x)\,dx
\right)^2 \\
& \quad 
+8C_2 
\sum_{1\le j<k\le n}
\left(
\int_{\textbf{R}^n} x_j^2 u_1(x)\,dx
\right)
\left(
\int_{\textbf{R}^n} x_k^2 u_1(x)\,dx
\right) 
+16C_2 
\sum_{1\le j<k\le n}
\left(
\int_{\textbf{R}^n} x_j x_k u_1(x)\,dx
\right)^2.
\end{align*}
Next, we treat $B_0^1(t,\xi)$. 
Similar calculations show 
\begin{align*}
\left\|
B_0^1(t,\xi) 
\right\|_{L^2(|\xi|\le1)}^2 
& \ge t^{-\frac{n}{4}-\frac{1}{2}}
\sum_{j=1}^n 
\left(
\int_{\textbf{R}^n} x_j u_0(x)\,dx
\right)^2
\int_{|\eta|\le1}
\eta_j^2
e^{-|\eta|^4} 
\cos^2(t^\frac{3}{4}|\eta|)
\,d\eta \\
& \ge \frac{1}{n}
t^{-\frac{n}{4}-\frac{1}{2}}
\sum_{j=1}^n 
\left(
\int_{\textbf{R}^n} x_j u_0(x)\,dx
\right)^2
\int_{|\eta|\le1}
|\eta|^2
e^{-|\eta|^4} 
\cos^2(t^\frac{3}{4}|\eta|)
\,d\eta \\
& \ge 2C_3
\sum_{j=1}^n 
\left(
\int_{\textbf{R}^n} x_j u_0(x)\,dx
\right)^2
t^{-\frac{n}{4}-\frac{1}{2}}, 
\qquad
t\gg1. 
\end{align*}
Therefore, we obtain 
\begin{align*}
& \left\|
B_1^2(t,\xi)+B_0^1(t,\xi) 
\right\|_{L^2(|\xi|\le1)}^2 
=\left\|
B_1^2(t,\xi) 
\right\|_{L^2(|\xi|\le1)}^2 
+\left\|
B_0^1(t,\xi) 
\right\|_{L^2(|\xi|\le1)}^2 \\
& \ge 2\left\{
C_1
\sum_{j=1}^n 
\left(
\int_{\textbf{R}^n} x_j^2 u_1(x)\,dx
\right)^2
+C_2\sum_{1\le j<k\le n}
\left(
\int_{\textbf{R}^n} x_j^2 u_1(x)\,dx
\right)
\left(
\int_{\textbf{R}^n} x_k^2 u_1(x)\,dx 
\right) 
\right. \\
& \qquad\qquad
\left. 
+2C_2 
\sum_{1\le j<k\le n}
\left(
\int_{\textbf{R}^n} x_j x_k u_1(x)\,dx
\right)^2
+C_3
\sum_{j=1}^n 
\left(
\int_{\textbf{R}^n} x_j u_0(x)\,dx
\right)^2
\right\}
t^{-\frac{n}{4}-\frac{1}{2}}
\end{align*}
for sufficiently large $t$. 
This implies the desired estimate. 
$\hfill \Box$
\\
\begin{remark}
{\rm 
From the above proof, we can obtain the following lower bound in the case of $n=1$: 
\begin{align*}
& \left\|
B_1^2(t,\xi)+B_0^1(t,\xi) 
\right\|_{L^2(|\xi|\le1)}^2 \\
& \ge 2C
\left\{
\left(
\int_{-\infty}^\infty x^2 u_1(x)\,dx
\right)^2
+
\left(
\int_{-\infty}^\infty x u_0(x)\,dx
\right)^2
\right\}
t^{-\frac{1}{4}-\frac{1}{2}},
\qquad
t\gg1,
\end{align*}
where
\[
C:=\frac{1}{16}
\int_0^1 \eta^2 e^{-\eta^4}\,d\eta.
\]
So, when $n=1$, the lower bound in Theorem~\ref{thm:4.5} becomes simpler. 
}
\end{remark}

Even in low regularity cases, we can define $C^k(t,\xi)$. 
In such cases, the corresponding solution and the profiles $C^k(t,\xi)$ decay slowly. 
However, the integrability of the solution in the high-frequency region is compensated by pulling out slowly decaying profiles. \\
\textbf{Proof of Theorem~\ref{thm:4.6}.} 
Let $j\in\textbf{N}_0$. 
By the Taylor theorem we see that 
\[
\begin{split}
&E_0(t,\xi)-\sum_{p=0}^j \mathcal{H}_0^{p}(t,\xi)
=\frac{1}{(j+1)!}
\frac{\partial^{j+1} H_0}{\partial b^{j+1}}(t,\xi,\theta|\xi|^{-6}) 
\cdot
\frac{1}{|\xi|^{6j+6}}
+\frac{1}{2}e^{\lambda_2 t},\\[2mm]
&E_1(t,\xi)-\sum_{p=1}^j \mathcal{H}_1^{p}(t,\xi)
=\frac{1}{(j+1)!}
\frac{\partial^{j+1} H_1}{\partial b^{j+1}}(t,\xi,\theta|\xi|^{-6}) 
\cdot
\frac{1}{|\xi|^{6j+6}}
+\frac{e^{\lambda_2 t}}{\sqrt{|\xi|^8-4|\xi|^2}},
\end{split}
\]
for some $0\le\theta\le1$. 
Now, we confirm that the function $\sqrt{1-4b}$ and its derivatives are all bounded for $0\le b\le1/8$. 
For $j\in\textbf{N}_0$, there exist constants $C>0$ and $c>0$ such that 
\[
\left|
\frac{\partial^j H_0}{\partial b^j}(t,\xi,b)
\right|
\le 
C e^{-\frac{ct}{|\xi|^2}}
\sum_{p=0}^j 
\left(
\frac{t}{|\xi|^2}
\right)^p, 
\qquad
\left|
\frac{\partial^j H_1}{\partial b^j}(t,\xi,b)
\right|
\le 
C \frac{e^{-\frac{ct}{|\xi|^2}}}{|\xi|^4}
\sum_{p=0}^j 
\left(
\frac{t}{|\xi|^2}
\right)^p, 
\]
for $t>0$, $\xi\in\textbf{R}^n$ and $0\le b\le1/8$. 
First, we consider the case of $k\in2\textbf{N}_0+1$, i.e., $k=2k'+1$ with $k'\in\textbf{N}_0$. 
It follows that 
\[
\begin{split}
&|\hat{u}(t,\xi)-C^k(t,\xi)| \\[2mm]
\le &
C \left(
\left|
\frac{\partial^{k'+1} H_0}{\partial b^{k'+1}}(t,\xi,\theta|\xi|^{-6})
\right|
\cdot 
\frac{1}{|\xi|^{6k'+6}} 
+\left|
\frac{\partial^{k'+1} H_1}{\partial b^{k'+1}}(t,\xi,\theta|\xi|^{-6})
\right|
\cdot 
\frac{1}{|\xi|^{6k'+2}}
\right)
|\widehat{u_0}(\xi)|  \\[2mm]
&+C\left|
\frac{\partial^{k'+1} H_1}{\partial b^{k'+1}}(t,\xi,\theta|\xi|^{-6})
\right|
\cdot 
\frac{1}{|\xi|^{6k'+6}} 
\cdot 
|\widehat{u_1}(\xi)|
+C e^{\lambda_2 t} (|\widehat{u_0}(\xi)|+|\widehat{u_1}(\xi)|) \\[2mm]
\le &
C \frac{e^{-\frac{ct}{|\xi|^2}}}{|\xi|^{3k+3}} 
\sum_{p=0}^{k'+1} 
\left(\frac{t}{|\xi|^2} \right)^p 
|\widehat{u_0}(\xi)|
+C \frac{e^{-\frac{ct}{|\xi|^2}}}{|\xi|^{3k+7}} 
\sum_{p=0}^{k'+1} \left(\frac{t}{|\xi|^2} \right)^p 
|\widehat{u_1}(\xi)| +C e^{\lambda_2 t} (|\widehat{u_0}(\xi)|+|\widehat{u_1}(\xi)|).
\end{split}
\]
Together with \eqref{7.2}, we obtain 
\[
\begin{split}
&\int_{|\xi| \ge \sqrt{2}} |\hat{u}(t,\xi)-C^k(t,\xi)|^2 d\xi \\[2mm]
\le& 
C\|u_0\|_{H^{l+1}}^2(1+t)^{-l-3k-4}
+C\|u_1\|_{H^l}^2(1+t)^{-l-3k-7}
+C \left(\|u_0\|_2^2+\|u_1\|_2^2 \right)\, e^{-\eta t} \\[2mm]
\le&
C \left(\|u_0\|_{H^{l+1}}^2+\|u_1\|_{H^l}^2 \right) (1+t)^{-l-3k-4},
\qquad
t\ge0.
\end{split}
\]
The case of $k\in2\textbf{N}_0$ is similar and so the details are reader. 
$\hfill \Box$
\\

\noindent
\textbf{Proof of Theorem~\ref{thm:4.7}.} 
It follows from Theorem~\ref{thm:4.4} and \eqref{4.6} with $k=0$ that 
\begin{align*}
& \left\|
\hat{u}(t,\cdot)
-A_1^0(t,\cdot)
-C^0(t,\cdot)
\right\|_2 \\
\le & \left\|
\hat{u}(t,\xi)
-A_1^0(t,\xi)
\right\|_{L^2(|\xi|\le1)} \\
& \qquad
+\left\|
\hat{u}(t,\xi)
-C^0(t,\xi)
\right\|_{L^2(|\xi|\ge\sqrt{2})}
+C(\|u_0\|_2+\|u_1\|_2+\|u_0\|_1+\|u_1\|_{1,1})e^{-\eta t} \\
\le & C(\|u_0\|_1+\|u_1\|_{1,1})(1+t)^{-\frac{n}{8}}
+C(\|u_0\|_{H^{l+1}}+\|u_1\|_{H^l})(1+t)^{-\frac{l+4}{2}},
\qquad
t\ge0.
\end{align*}
We also see that 
\[
\|C^0(t,\cdot)\|_2
\le C\|u_0\|_{H^{l+1}}
\cdot
t^{-\frac{l+1}{2}},
\qquad t>0. 
\]
With a slight modification of \eqref{7.4}, \eqref{7.5} and \eqref{7.6} (or see \cite{II} directly), we have 
\begin{align*}
\|A_1^0(t,\cdot)\|_2
\le 
\begin{cases}
C\|u_0\|_1
\sqrt{t}, 
& n=1, \\[7pt]
C\|u_0\|_1
\sqrt{\log t}, 
& n=2, \\[7pt]
C\|u_0\|_1
t^{-\frac{n}{8}+\frac{1}{4}}, 
& n\ge3, 
\end{cases}
\qquad t>0.  
\end{align*}
(i) It suffices to consider the case of $n\ge3$. 
Then, one has 
\begin{align*}
\left\|
\hat{u}(t,\cdot)
-C^0(t,\cdot)
\right\|_2 
& \le C(\|u_0\|_1+\|u_1\|_{1,1})(1+t)^{-\frac{n}{8}+\frac{1}{4}}
+C(\|u_0\|_{H^{l+1}}+\|u_1\|_{H^l})t^{-\frac{l+4}{2}},
\qquad
t>0.
\end{align*}
When 
\[
-\frac{n}{8}+\frac{1}{4}
<-\frac{l+1}{2},
\qquad
\mbox{i.e.,}
\qquad
l<\frac{n}{4}-\frac{3}{2},
\]
the function $C^0(t,\xi)$ is the leading term. 
If 
\[
-\frac{n}{8}+\frac{1}{4}
\le-\frac{l+4}{2},
\qquad
\mbox{i.e.,}
\qquad
l\le\frac{n}{4}-\frac{9}{2},
\] 
then it holds
\begin{align*}
\left\|
\hat{u}(t,\cdot)
-C^0(t,\cdot)
\right\|_2 
& \le CI^{l,0}(u_0,u_1) t^{-\frac{l+4}{2}},
\qquad
t>0.
\end{align*}
On the other hand, if 
\[
-\frac{l+4}{2}< -\frac{n}{8}+\frac{1}{4},
\qquad
\mbox{i.e.,}
\qquad
\frac{n}{4}-\frac{9}{2}<l,
\] 
then we have 
\begin{align*}
\left\|
\hat{u}(t,\cdot)
-C^0(t,\cdot)
\right\|_2 
& \le CI^{l,0}(u_0,u_1) 
t^{-\frac{n}{8}+\frac{1}{4}},
\qquad
t>0.
\end{align*}
So, under the condition $n\ge3$, we obtain 
\begin{align*}
\left\|
\hat{u}(t,\xi)
-e^{-\frac{t}{|\xi|^2}}
\widehat{u_0}
\right\|_2 
& \le 
\begin{cases}
CI^{l,0}(u_0,u_1) 
t^{-\frac{l+4}{2}},
& l\le n/4-9/2, \\[7pt]
CI^{l,0}(u_0,u_1) 
t^{-\frac{n}{8}+\frac{1}{4}},
& n/4-9/2<l<n/4-3/2, 
\end{cases}
\end{align*}
for $t\ge1$. 
This gives statement~(i). \\
(ii) Note that 
\[
-\frac{l+4}{2}<-\frac{n}{8}
\qquad
\mbox{if}
\qquad
l=\frac{n}{4}-\frac{3}{2}.
\]
Hence, statement~(ii) easily follows. \\
(iii) We have 
\begin{align*}
\left\|
\hat{u}(t,\cdot)
-A_1^0(t,\cdot)
\right\|_2 
& \le C(\|u_0\|_1+\|u_1\|_{1,1})(1+t)^{-\frac{n}{8}}
+C(\|u_0\|_{H^{l+1}}+\|u_1\|_{H^l})t^{-\frac{l+1}{2}},
\qquad
t>0.
\end{align*}
In the case of $n\ge3$, if  
\[
-\frac{l+1}{2}
<-\frac{n}{8}
+\frac{1}{4}, 
\qquad
\mbox{i.e.,}
\qquad
\frac{n}{4}-\frac{3}{2}<l,
\]
the function $A_1^0(t,\xi)$ is the leading term. 
Furthermore, if 
\[
-\frac{n}{8}
\le -\frac{l+1}{2},
\qquad
\mbox{i.e.,}
\qquad
l\le \frac{n}{4}-1,
\]
then we have 
\begin{align*}
\left\|
\hat{u}(t,\cdot)
-A_1^0(t,\cdot)
\right\|_2 
& \le C I^{l,0}(u_0,u_1) t^{-\frac{l+1}{2}},
\qquad
t\ge1.
\end{align*}
Conversely, when 
\[
-\frac{l+1}{2}
<-\frac{n}{8},
\qquad
\mbox{i.e.,}
\qquad
\frac{n}{4}-1<l,
\]
we have 
\begin{align*}
\left\|
\hat{u}(t,\cdot)
-A_1^0(t,\cdot)
\right\|_2 
& \le C I^{l,0}(u_0,u_1) t^{-\frac{n}{8}},
\qquad
t\ge1.
\end{align*}
Thus, when $n\ge3$, we obtain 
\begin{align*}
& \left\|
\hat{u}(t,\xi)
-\left(
\int_{\textbf{R}^n} u_1(x)\,dx
\right)
e^{-\frac{1}{2}t|\xi|^4}
\frac{\sin(t|\xi|)}{|\xi|}
\right\|_2 
\le 
\begin{cases}
C I^{l,0}(u_0,u_1) 
t^{-\frac{l+1}{2}}, 
& n/4-3/2<l\le n/4-1, \\[7pt]
C I^{l,0}(u_0,u_1) 
t^{-\frac{n}{8}}, 
& n/4-1<l,
\end{cases}
\end{align*}
for $t\ge1$. 
In the case of $n=1,2$, we can easily see that, for $l\ge0$, 
\begin{align*}
\left\|
\hat{u}(t,\xi)
-\left(
\int_{\textbf{R}^n} u_1(x)\,dx
\right)
e^{-\frac{1}{2}t|\xi|^4}
\frac{\sin(t|\xi|)}{|\xi|}
\right\|_2 
\le 
C I^{l,0}(u_0,u_1) 
t^{-\frac{n}{8}},
\qquad
t\ge1. 
\end{align*}
Thus, we obtain (iii) and the proof is now complete. 
$\hfill \Box$
\\

Corollary~\ref{cor:4.1} is a direct result of Theorem~\ref{thm:4.3} with Lemmas~\ref{lem:7.3}, \ref{lem:7.4} and \ref{lem:7.5}. 
\\

%%%%%%%%%%%%%%%%%%%%%%%%%%%%%%%%%%%%%%%%%%%%%%%%%%%%%%%
%%%%%%%%%%%%%%
\section{Application}
\label{sec:6}
%%%%%%%%%%%%%%
\begin{theorem}
\label{thm:6.1}
Let $n\ge3$ and $u$ be the solution to \eqref{1.1} with $u_0\in H^1(\textbf{R}^n)\cap L^1(\textbf{R}^n)$, $u_1\in L^2(\textbf{R}^n)\cap L^{1,1}(\textbf{R}^n)$. 
If $(n-18)/12\le k\in\textbf{N}_0$, then it holds
\begin{align*}
 C_1\left|
\int_{\textbf{R}^n} u_1(x)\,dx
\right|
t^{-\frac{n}{8}+\frac{1}{4}}
& \le \|\hat{u}(t,\cdot)-C^k(t,\cdot)\|_2
\le C_2 I^{0,0}(u_0,u_1) 
t^{-\frac{n}{8}+\frac{1}{4}}, 
\qquad
t\gg1.
\end{align*}
Here, $C_1>0$ and $C_2>0$  are constants independent of $t$ and the initial data. 
\end{theorem}
\textbf{Proof.}
It suffices to show the upper bound. 
It follows from \eqref{4.4} with $\gamma=1$, \eqref{7.6} and \eqref{4.6} with $l=0$ that 
\[
\begin{split}
&\|\hat{u}(t,\cdot)-C^k(t,\cdot)\|_2 \\[2mm]
\le&
\|\hat{u}(t,\xi)\|_{L^2(|\xi| \le 1)}
+\|\hat{u}(t,\xi)-C^k(t,\xi)\|_{L^2(|\xi|\ge\sqrt{2})}
+C(\|u_0\|_2+\|u_1\|_2)e^{-\eta t} \\[2mm]
\le&
C(\|u_0\|_1+\|u_1\|_{1,1})t^{-\frac{n}{8}+\frac{1}{4}}
+C(\|u_0\|_{H^1}+\|u_1\|_2)t^{-\frac{3k+4}{2}}
+C(\|u_0\|_2+\|u_1\|_2)e^{-\eta t} \\[2mm]
\le&
CI^{0,0}(u_0,u_1) 
t^{-\frac{n}{8}+\frac{1}{4}}, \qquad t \gg 1,
\end{split}
\]
for some $\eta>0$. 
The condition on $k$ is used to derive the last inequality. 
$\hfill \Box$

\begin{remark}
{\rm 
In addition to assumptions in Theorem~\ref{thm:6.1}, if  $(u_0,u_1)\in H^{l+1}(\textbf{R}^n)\times H^l(\textbf{R}^n)$ with $l\ge l^*=n/4-3/2$ and $l\ge0$, the following estimate is shown (recall Corollary~\ref{cor:4.1} and Theorem~\ref{thm:4.3}, and see also \cite[Theorem 4.1]{II} with Remark~\ref{rem:4.3}): 
\[
C_1 |P_1|
t^{-\frac{n}{8}+\frac{1}{4}}
\le \|u(t,\cdot)\|_2
\le C_2 I^{l,0}(u_0,u_1) t^{-\frac{n}{8}+\frac{1}{4}},
\quad
\mbox{where}
\quad
P_1:=\int_{\textbf{R}^n} u_1(x)\,dx,
\]
for sufficiently large $t$. 
In the low dimensional cases $1\le n\le6$, the diffusive structure is strong even for the energy solution. 
However, we have to impose additional regularity condition $l\ge l^*$ in the case of $n\ge7$ to expect to observe the same estimate. 
In this sense, Theorem~\ref{thm:6.1} is essential in the case of $n\ge7$,  and it says that similar effects can be obtained by pulling out the slowly decay profiles $C^k(t,\xi)$ from the Fourier transformed solution instead of taking smooth initial data. 
We insist that, for higher dimensional cases, $C^k(t,\xi)$ with 
\[
k^*:=\frac{n-18}{12}\le k\in\textbf{N}_0
\]
is a non-diffusive part of the solution. 
The value $k^*$ indicates the critical point of whether the difference $\hat{u}(t,\xi)-C^k(t,\xi)$ with $k\ge k^*$ possesses the diffusion wave property. 
This threshold $k^*$ is one of novelties in this study. 
For example, Theorem~\ref{thm:6.1} gives 
\begin{align*}
 C_1|P_1|
t^{-\frac{n}{8}+\frac{1}{4}}
& \le \left\|
\hat{u}(t,\xi)-e^{-\frac{t}{|\xi|^2}}\widehat{u_0}
\right\|_2
\le C_2 I^{0,0}(u_0,u_1) 
t^{-\frac{n}{8}+\frac{1}{4}}, 
\qquad
7\le n\le 18,
\end{align*}
\begin{align*}
 C_1|P_1|
 t^{-\frac{n}{8}+\frac{1}{4}}
& \le \left\|
\hat{u}(t,\xi)
-e^{-\frac{t}{|\xi|^2}}\widehat{u_0}
-\frac{e^{-\frac{t}{|\xi|^2}}}{|\xi|^4}\widehat{u_1}
\right\|_2
\le C_2 I^{0,0}(u_0,u_1) 
t^{-\frac{n}{8}+\frac{1}{4}}, 
\qquad
19\le n\le 30,
\end{align*}
for sufficiently large $t$. 
}
\end{remark}

Similarly, we obtain the following theorems. 
\begin{theorem}
Let $n\ge3$ and $u$ be the solution to \eqref{1.1} with $u_0\in H^1(\textbf{R}^n)\cap L^1(\textbf{R}^n)$, $u_1\in L^2(\textbf{R}^n)\cap L^{1,1}(\textbf{R}^n)$. 
If $(n-16)/12\le k\in\textbf{N}_0$, then it holds
\begin{align*}
& C_1 \sqrt{
\left(
\int_{\textbf{R}^n} u_0(x)\,dx
\right)^2
+\sum_{j=1}^n
\left(
\int_{\textbf{R}^n} x_j u_1(x)\,dx
\right)^2
}t^{-\frac{n}{8}} \\
& \qquad \le \left\|
\hat{u}(t,\xi)
-\left(
\int_{\textbf{R}^n} u_1(x)\,dx
\right)
e^{-\frac{1}{2}t|\xi|^4}
\frac{\sin(t|\xi|)}{|\xi|}
-C^k(t,\xi)
\right\|_2
\le C_2 I^{0,0}(u_0,u_1) t^{-\frac{n}{8}}
\end{align*}
for sufficiently large $t$. 
Here, $C_1>0$ and $C_2>0$  are constants independent of $t$ and the initial data. 
\end{theorem}

\begin{theorem}
Let $n\ge3$ and $u$ be the solution to \eqref{1.1} with $u_0\in H^1(\textbf{R}^n)\cap L^{1,1}(\textbf{R}^n)$, $u_1\in L^2(\textbf{R}^n)\cap L^{1,2}(\textbf{R}^n)$. 
If $(n-14)/12\le k\in\textbf{N}_0$, then it holds 
\begin{align*}
& \left\{
C_1
\sum_{j=1}^n 
\left(
\int_{\textbf{R}^n} x_j^2 u_1(x)\,dx
\right)^2
+C_2\sum_{1\le j<k\le n}
\left(
\int_{\textbf{R}^n} x_j^2 u_1(x)\,dx
\right)
\left(
\int_{\textbf{R}^n} x_k^2 u_1(x)\,dx 
\right) 
\right. \\
& \qquad\qquad
\left. 
+2C_2 
\sum_{1\le j<k\le n}
\left(
\int_{\textbf{R}^n} x_j x_k u_1(x)\,dx
\right)^2
+C_3
\sum_{j=1}^n 
\left(
\int_{\textbf{R}^n} x_j u_0(x)\,dx
\right)^2
\right\}^\frac{1}{2}
t^{-\frac{n}{8}-\frac{1}{4}} \\
\le & \left\|
\hat{u}(t,\xi)
-\left(
\int_{\textbf{R}^n} u_1(x)\,dx
\right)
e^{-\frac{1}{2}t|\xi|^4}
\frac{\sin(t|\xi|)}{|\xi|} 
\right.\\
& \quad \left.
+i \sum_{j=1}^n
\left(
\int_{\textbf{R}^n} x_j u_1(x)\,dx
\right)
e^{-\frac{1}{2}t|\xi|^4}
\frac{\sin(t|\xi|)}{|\xi|}
\xi_j
-\left(
\int_{\textbf{R}^n} u_0(x)\,dx
\right)
e^{-\frac{1}{2}t|\xi|^4}
\cos(t|\xi|)
-C^k(t,\xi)
\right\|_2 \\
\le & C I^{0,1}(u_0,u_1) 
t^{-\frac{n}{8}-\frac{1}{4}}
\end{align*}
for sufficiently large $t$. 
Here, $C>0$ is a constant independent of $t$ and the initial data and $C_j$ $(j=1,2,3)$ are defined in Theorem~\ref{thm:4.5}. 
\end{theorem}

%%%%%%%%%%%%%%
\section{Appendix}
\label{sec:7}
%%%%%%%%%%%%%%
\begin{lemma}
Let $n\ge1$, $m\ge0$ and $c>0$. 
Then, there exists a constant $C>0$ such that 
\begin{align}
\label{7.1}
\big\||\xi|^m e^{-ct|\xi|^4}\big\|_{L^2(|\xi|\le1)}
\le C(1+t)^{-\frac{n}{8}-\frac{m}{4}},
\qquad 
t\ge0,
\end{align}
%%%
\begin{align}
\label{7.2}
\sup_{|\xi|\ge1}
\big(|\xi|^{-m} e^{-\frac{ct}{|\xi|^2}}\big)
\le C(1+t)^{-\frac{m}{2}},
\qquad 
t\ge0.
\end{align}
\end{lemma}

\begin{lemma}
{\rm (\cite{IM, Mi4})}
\label{lem:7.2}
Let $n\ge1$ and $f\in L^{1,\gamma}(\textbf{R}^n)$ with $\gamma\ge0$. 
Then, it holds
\begin{align}
\label{7.3}
\left|
\hat{f}(\xi)
-\sum_{|\alpha|\le[\gamma]}
M_\alpha(f)(i\xi)^\alpha
\right|
\le 
C |\xi|^\gamma \|f\|_{1,\gamma}, 
\qquad
\xi\in\textbf{R}^n.
\end{align}
Here, $C>0$ is a constant independent of $\xi$ and $f$. 
\end{lemma}

In \cite{II}, the optimal estimate for the integral 
\[
\int_{\textbf{R}^n} 
e^{-t|\xi|^4}
\frac{\sin^2(t|\xi|)}{|\xi|^2}
\,d\xi
\]
for all $n\ge1$ based on \cite{IO}. 
Together with the following estimate
\[
\int_{|\xi|\ge1}
e^{-t|\xi|^4}
\frac{\sin^2(t|\xi|)}{|\xi|^2}
\,d\xi
\le \int_{|\xi|\ge1}
e^{-t|\xi|^4}
\,d\xi
\le \left(
\int_{|\xi|\ge1}
e^{-\frac{|\xi|^4}{2}}
\,d\xi
\right)
e^{-\frac{t}{2}},
\qquad
t\ge1,
\]
we can prove Lemmas~\ref{7.3}, \ref{7.4} and \ref{7.5} immediately.
However, in the latter part of this section, we give direct proofs that show the same estimates also hold even if we replace the integral domain from $\textbf{R}^n$ to $\{\xi\in\textbf{R}^n:|\xi|\le1\}$. 
\begin{lemma}
\label{lem:7.3}
Let $n=1$. 
There exist constants $0<C_1^1<C_2^1<\infty$ such that 
\begin{align}
\label{7.4}
C_1^1 t
\le \int_{-1}^1
e^{-t\xi^4}
\frac{\sin^2(t|\xi|)}{\xi^2}
\,d\xi
\le C_2^1 t, 
\qquad
t\gg1.
\end{align}
\end{lemma}
\textbf{Proof.} 
First, we see that 
\[
\int_{-1}^1
e^{-t\xi^4}
\frac{\sin^2(t|\xi|)}{\xi^2}
\,d\xi
=2t^\frac{1}{4}
\int_0^{t^\frac{1}{4}}
e^{-\eta^4}
\frac{\sin^2(t^\frac{3}{4}\eta)}{\eta^2}
\,d\eta.
\]
Since 
\[
\frac{\theta}{2}
\le \sin\theta
\le\theta
\qquad
\mbox{for}
\qquad
0\le\theta\le1,
\]
we have 
\[
\frac{t^\frac{3}{4}\eta}{2}
\le \sin(t^\frac{3}{4}\eta)
\le t^\frac{3}{4}\eta
\qquad
\mbox{for}
\qquad
0\le
\eta
\le t^{-\frac{3}{4}}.
\]
We divide the integral into two parts:
\begin{align*}
\int_0^{t^{\frac{1}{4}}}
e^{-\eta^4}
\frac{\sin^2(t^\frac{3}{4}\eta)}{\eta^2}
\,d\eta
=
\left(
\int_0^{t^{-\frac{3}{4}}}
+
\int_{t^{-\frac{3}{4}}}^{t^{\frac{1}{4}}}
\right)
e^{-\eta^4}
\frac{\sin^2(t^\frac{3}{4}\eta)}{\eta^2}
\,d\eta
=:I_1+I_2, 
\qquad
t\ge1.
\end{align*}
It follows that 
\begin{align*}
I_1\le t^{\frac{3}{2}}
\int_0^{t^{-\frac{3}{4}}}
e^{-\eta^4}
\,d\eta
\le  t^{\frac{3}{2}}
\cdot
t^{-\frac{3}{4}}
=t^{\frac{3}{4}}, 
\qquad
t\ge1,
\end{align*}
\begin{align*}
I_1
& \ge \frac{t^{\frac{3}{2}}}{4}
\int_0^{t^{-\frac{3}{4}}}
e^{-\eta^4}
\,d\eta
\ge \frac{t^{\frac{3}{2}}}{4}
\int_0^{t^{-\frac{3}{4}}}
e^{-\eta}
\,d\eta
\ge \frac{t^{\frac{3}{2}}}{4}
\int_{t^{-\frac{3}{4}}/2}^{t^{-\frac{3}{4}}}
e^{-\eta} 
\,d\eta \\
& \ge \frac{t^{\frac{3}{2}}}{4}
\cdot 
\frac{t^{-\frac{3}{4}}}{2}
\exp\left(-\frac{1}{t^{\frac{3}{4}}}\right)
\ge \frac{t^{\frac{3}{4}}}{8e},
\qquad
t\ge1.
\end{align*}
Next, we see that 
\begin{align*}
I_2
& \le \int_{t^{-\frac{3}{4}}}^{t^{\frac{1}{4}}}
\frac{e^{-\eta^4}}{\eta^2}
\,d\eta
=\int_{t^{-\frac{3}{4}}}^{t^{\frac{1}{4}}}
\left(-\frac{1}{\eta}\right)'
e^{-\eta^4}
\,d\eta \\
& \le \left[
-\frac{e^{-\eta^4}}{\eta}
\right]_{t^{-\frac{3}{4}}}^{t^{\frac{1}{4}}}
-4\int_{t^{-\frac{3}{4}}}^{t^{\frac{1}{4}}}
\eta^2 
e^{-\eta^4}
\,d\eta
\le t^{\frac{3}{4}}
\exp\left(-\frac{1}{t^3}\right)
\le t^{\frac{3}{4}}, 
\qquad 
t\ge1.
\end{align*}
Since 
\[
\sin\theta\ge\frac{1}{2}
\qquad
\mbox{for}
\qquad
1\le\theta\le\frac{\pi}{2},
\]
we have 
\[
\sin(t^{\frac{3}{4}}\eta)\ge\frac{1}{2}
\qquad
\mbox{for}
\qquad
t^{-\frac{3}{4}} \le \eta\le\frac{\pi}{2}t^{-\frac{3}{4}}.
\]
It follows that 
\begin{align*}
I_2
& \ge \frac{1}{4}
\int_{t^{-\frac{3}{4}}}^{\frac{\pi}{2}t^{-\frac{3}{4}}}
\frac{e^{-\eta^4}}{\eta^2}
\,d\eta
=\frac{1}{4}
\left[
-\frac{e^{-\eta^4}}{\eta}
\right]_{t^{-\frac{3}{4}}}^{\frac{\pi}{2}t^{-\frac{3}{4}}}
-\int_{t^{-\frac{3}{4}}}^{\frac{\pi}{2}t^{-\frac{3}{4}}}
\eta^2 e^{-\eta^4}
\,d\eta \\
& \ge 
\frac{1}{4}
t^{\frac{3}{4}}
\left(
\exp\left(-\frac{1}{t^3}\right)
-\frac{2}{\pi}
\exp\left(-\frac{\pi^4}{16 t^3}\right)
\right)
-\int_0^2 \eta^2 e^{-\eta^4}
\,d\eta \\
& =
\frac{1}{4}
t^{\frac{3}{4}}
\exp\left(-\frac{1}{t^3}\right)
\left(
1
-\frac{2}{\pi}
\exp\left(-\frac{c}{t^3}\right)
\right)
-C \\
& \ge 
\frac{1}{4e}
\left(
1-\frac{2}{\pi}
\right)
t^{\frac{3}{4}}
-C, 
\qquad
t\ge1, 
\end{align*}
where 
\[
c:=\frac{\pi^4}{16}-1>0, 
\qquad
C:=\int_0^2 \eta^2 e^{-\eta^4}
\,d\eta.
\]
Thus, one has 
\[
I_2
\ge \frac{1}{8e}
\left(
1-\frac{2}{\pi}
\right)
t^{\frac{3}{4}},
\qquad
t\gg1.
\]
The proof is now complete. 
$\hfill \Box$

\begin{lemma}
\label{lem:7.4}
Let $n=2$. 
There exist constants $0<C_1^2<C_2^2<\infty$ such that 
\begin{align}
\label{7.5}
C_1^2 \log t
\le \int_{|\xi|\le1}
e^{-t|\xi|^4}
\frac{\sin^2(t|\xi|)}{|\xi|^2}
\,d\xi
\le C_2^2 \log t, 
\qquad
t\gg1.
\end{align}
\end{lemma}
\textbf{Proof.} 
Similarly, we divide the integral as follows:
\begin{align*}
& \int_{|\xi|\le1}
e^{-t|\xi|^4}
\frac{\sin^2(t|\xi|)}{|\xi|^2}
\,d\xi
=\int_{|\eta|\le t^\frac{1}{4}}
e^{-|\eta|^4}
\frac{\sin^2(t^\frac{3}{4}|\eta|)}{|\eta|^2}
\,d\eta \\
= & \left(
\int_{|\eta|\le t^{-\frac{3}{4}}}
+\int_{t^{-\frac{3}{4}}\le |\eta|\le t^\frac{1}{4}}
\right)
e^{-|\eta|^4}
\frac{\sin^2(t^\frac{3}{4}|\eta|)}{|\eta|^2}
\,d\eta
=:I_3+I_4, 
\qquad
t\ge1.
\end{align*}
It follows that 
\begin{align*}
I_3
\le t^{\frac{3}{2}}
\int_{|\eta|\le t^{-\frac{3}{4}}}
e^{-|\eta|^4}
\,d\eta
\le t^{\frac{3}{2}}
\int_{|\eta|\le t^{-\frac{3}{4}}}
\,d\eta
=\pi,
\qquad
t\ge1,
\end{align*}
\begin{align*}
I_3
& \ge \frac{t^{\frac{3}{2}}}{4}
\int_{|\eta|\le t^{-\frac{3}{4}}}
e^{-|\eta|^4}
\,d\eta
=\frac{\pi}{2}
t^{\frac{3}{2}}
\int_0^{t^{-\frac{3}{4}}}
r e^{-r^4}
\,dr
\ge 
\frac{\pi}{2}
t^{\frac{3}{2}}
\int_0^{t^{-\frac{3}{4}}}
r^3 e^{-r^4}
\,dr \\
& =\frac{\pi}{2}
t^{\frac{3}{2}}
\left[
-\frac{1}{4}e^{-r^4}
\right]_0^{t^{-\frac{3}{4}}}
=\frac{\pi}{8}
t^{\frac{3}{2}}
\left(
1-\exp\left(-\frac{1}{t^3}\right)
\right)
=\frac{\pi}{8}
t^{-\frac{3}{2}}
+O(t^{-\frac{9}{2}}),
\qquad 
t\gg1,
\end{align*}
since 
\begin{align*}
1-\exp\left(-\frac{1}{t^3}\right)
=-\sum_{k=1}^\infty
\frac{1}{k!}
\left(
-\frac{1}{t^3}
\right)^k
=\frac{1}{t^3}
-\sum_{k=2}^\infty
\frac{1}{k!}
\left(
-\frac{1}{t^3}
\right)^k
\ge \frac{1}{t^3}
-\left(
\sum_{k=2}^\infty
\frac{1}{k!}
\right)
\frac{1}{t^6}.
\end{align*}
Next, we see that 
\begin{align*}
I_4
& \le \int_{t^{-\frac{3}{4}}\le |\eta|\le t^\frac{1}{4}}
\frac{e^{-|\eta|^4}}{|\eta|^2}
\,d\eta
=2\pi 
\int_{t^{-\frac{3}{4}}}^{t^\frac{1}{4}}
\frac{e^{-r^4}}{r}
\,dr
=2\pi 
\int_{t^{-\frac{3}{4}}}^{t^\frac{1}{4}}
(\log r)'
e^{-r^4}
\,dr \\
& = 2\pi
\left[
(\log r)e^{-r^4}
\right]_{t^{-\frac{3}{4}}}^{t^\frac{1}{4}}
+8\pi 
\int_{t^{-\frac{3}{4}}}^{t^\frac{1}{4}}
(\log r)r^3
e^{-r^4}
\,dr \\
& =2\pi\left\{
\frac{1}{4}
(\log t)
\exp\left(-\frac{1}{t}\right)
+\frac{3}{4}
(\log t)
\exp\left(-\frac{1}{t^3}\right)
\right\}
+8\pi 
\int_{t^{-\frac{3}{4}}}^{t^\frac{1}{4}}
(\log r)r^3
e^{-r^4}
\,dr, \\
& \le 2\pi\log t
+8\pi C_0, 
\qquad
t\ge1,
\end{align*}
where
\[
C_0
:=\int_0^\infty
|\log r| r^3
e^{-r^4}
\,dr
<+\infty.
\]
Note that the function $e^{-x^4}/x^2$ $(x>0)$ is monotone decreasing and 
\[
\sin^2(t^{\frac{3}{4}}|\eta|)
\ge \frac{1}{2} 
\qquad
\mbox{for}
\qquad
\left(
\frac{\pi}{4}
+j\pi
\right)
t^{-\frac{3}{4}}
\le |\eta|
\le \left(
\frac{3\pi}{4}
+j\pi
\right)
t^{-\frac{3}{4}}
\quad
(j\in\textbf{N}).
\]
The largest number $j\in\textbf{N}$ satisfying $(3\pi/4+j\pi)t^{-3/4}\le t^{1/4}$ is $[t/\pi-3/4]$. 
Thus, one has
\begin{align*}
I_4
& \ge \frac{1}{2}
\sum_{j=1}^{[t/\pi-3/4]}
\int_{(\frac{\pi}{4}+j\pi)t^{-\frac{3}{4}}\le |\eta|\le (\frac{3\pi}{4}+j\pi)t^{-\frac{3}{4}}}
\frac{e^{-|\eta|^4}}{|\eta|^2}
\,d\eta \\
& \ge \frac{1}{4}
\int_{(\frac{5\pi}{4}+j\pi)t^{-\frac{3}{4}}\le |\eta|\le (\frac{5\pi}{4}+[t/\pi-3/4]\pi)t^{-\frac{3}{4}}}
\frac{e^{-|\eta|^4}}{|\eta|^2}
\,d\eta \\
& \ge \frac{\pi}{2}
\int_{\frac{5\pi}{4}t^{-\frac{3}{4}}}^{(\frac{5\pi}{4}+[t/\pi-3/4]\pi)t^{-\frac{3}{4}}}
\frac{e^{-r^4}}{r}
\,dr \\
& =\frac{\pi}{2}
\left[
(\log r)e^{-r^4}
\right]_{\frac{5\pi}{4}t^{-\frac{3}{4}}}^{(\frac{5\pi}{4}+[t/\pi-3/4]\pi)t^{-\frac{3}{4}}}
+2\pi
\int_{\frac{5\pi}{4}t^{-\frac{3}{4}}}^{(\frac{5\pi}{4}+[t/\pi-3/4]\pi)t^{-\frac{3}{4}}}
(\log r)r^3 e^{-r^4}
\,dr \\
& =\frac{\pi}{2}
\left(
\frac{3}{4}\log t
-\log\left(
\frac{5\pi}{4}
\right)
\right)
\exp\left(
-\left(\frac{5\pi}{4}\right)^4
\frac{1}{t^3}
\right) \\
& \qquad
-\frac{\pi}{2}
\left(
\frac{3}{4}\log t 
-\log\left(
\frac{5\pi}{4}+\left[\frac{t}{\pi}-\frac{3}{4}\right]\pi
\right)
\right)
\exp\left(
-\left(
\frac{5\pi}{4}+\left[\frac{t}{\pi}-\frac{3}{4}\right]\pi
\right)^4
\frac{1}{t^3}
\right)
-2\pi C_0 \\
& =\frac{3\pi}{8}
(\log t)
\exp\left(
-\left(\frac{5\pi}{4}\right)^4
\frac{1}{t^3}
\right)
\left(1-\exp\left(-\frac{\alpha}{t^3}\right)\right)
-\frac{\pi}{2}
\log\left(
\frac{5\pi}{4}
\right)
\exp\left(
-\left(\frac{5\pi}{4}\right)^4
\frac{1}{t^3}
\right) \\
& \qquad
+\frac{\pi}{2}
\log\left(
\frac{5\pi}{4}+\left[\frac{t}{\pi}-\frac{3}{4}\right]\pi
\right)
\exp\left(
-\left(
\frac{5\pi}{4}+\left[\frac{t}{\pi}-\frac{3}{4}\right]\pi
\right)^4
\frac{1}{t^3}
\right)
-2\pi C_0, 
\qquad
t\gg1,
\end{align*}
where 
\begin{align*}
\alpha
:=\left(
\frac{5\pi}{4}+\left[\frac{t}{\pi}-\frac{3}{4}\right]\pi
\right)^4
-\left(
\frac{5\pi}{4}
\right)^4.
\end{align*}
Since 
\begin{align*}
\frac{\alpha}{t^3}
\ge \left[\frac{t}{\pi}-\frac{3}{4}\right]^4 \frac{\pi^4}{t^3}, 
\end{align*}
there exists a constant $c'>0$ such that 
\[
\exp\left(-\frac{\alpha}{t^3}\right)
\le  e^{-c't}, 
\qquad
t\gg1.
\]
Hence, we obtain 
\begin{align*}
I_4
& \ge \frac{3\pi}{16}
\exp\left(
-\left(\frac{5\pi}{4}\right)^4
\right)
(\log t)
-\frac{\pi}{2}
\log\left(
\frac{5\pi}{4}
\right) 
-2\pi C_0, 
\qquad
t\gg1.
\end{align*}
Therefore, we complete the proof. 
$\hfill \Box$

\begin{lemma}
\label{lem:7.5}
Let $n\ge3$. 
There exist constants $0<C_1^n<C_2^n<\infty$ such that 
\begin{align}
\label{7.6}
C_1^n t^{-\frac{n}{4}+\frac{1}{2}}
\le \int_{|\xi|\le1}
e^{-t|\xi|^4}
\frac{\sin^2(t|\xi|)}{|\xi|^2}
\,d\xi
\le C_2^n t^{-\frac{n}{4}+\frac{1}{2}}, 
\qquad
t\gg1.
\end{align}
\end{lemma}
\textbf{Proof.} 
It follows that 
\begin{align*}
\int_{|\xi|\le1}
e^{-t|\xi|^4}
\frac{\sin^2(t|\xi|)}{|\xi|^2}
\,d\xi
=t^{-\frac{n}{4}+\frac{1}{2}}
\int_{|\eta|\le t^\frac{1}{4}}
e^{-|\eta|^4}
\frac{\sin^2(t^\frac{3}{4}|\eta|)}{|\eta|^2}
\,d\eta. 
\end{align*}
We can easily see that 
\[
\int_{|\eta|\le t^\frac{1}{4}}
e^{-|\eta|^4}
\frac{\sin^2(t^\frac{3}{4}|\eta|)}{|\eta|^2}
\,d\eta
\le \int_{\textbf{R}^n}
\frac{e^{-|\eta|^4}}{|\eta|^2}
\,d\eta
=\omega_n
\int_0^\infty
x^{n-3}
e^{-x^4}
\,dx
<+\infty, 
\qquad
t>0.
\]
On the other hand, one has 
\begin{align*}
& \int_{|\eta|\le t^\frac{1}{4}}
e^{-|\eta|^4}
\frac{\sin^2(t^\frac{3}{4}|\eta|)}{|\eta|^2}
\,d\eta
\ge \int_{|\eta|\le1}
e^{-|\eta|^4}
\frac{\sin^2(t^\frac{3}{4}|\eta|)}{|\eta|^2}
\,d\eta \\
= & \frac{1}{2} 
\int_{|\eta|\le1}
\frac{e^{-|\eta|^4}}{|\eta|^2}
\,d\eta
-\frac{1}{2}
\int_{|\eta|\le1}
\frac{e^{-|\eta|^4}}{|\eta|^2}
\cos(2t^\frac{3}{4}|\eta|)
\,d\eta, 
\qquad
t\ge1.
\end{align*}
The Riemann-Lebesgue lemma gives the lower bound. 
$\hfill \Box$
\\

\noindent 
\textbf{Acknowledgement.} 
The work of the second author was supported in part by Grant-in-Aid for Scientific Research (C) 15K04958 of JSPS. 

%%%%%%%%%%%%%%%%%%%%%%%%%%%%%
\bibliographystyle{amsplain}

\end{document}